\documentclass{amsart}
\usepackage{titlesec}
\usepackage{amsfonts}
\usepackage{hhline}
\usepackage{amsmath}
\usepackage{enumerate}
\usepackage[shortlabels]{enumitem}
\usepackage[font={scriptsize}]{caption}
\usepackage{amsthm}
\usepackage{mathrsfs}
\usepackage{caption}
\usepackage{tikz}
\usepackage[all,cmtip]{xy}
\usepackage{amssymb}
\usepackage{changepage}
\usepackage{rotating}
\usepackage{array}
\usepackage{layout}
\def\it{\textit}

\def\mr{\mathbb{R}}
\def\mc{\mathbb{C}}
\def\mq{\mathbb{Q}}
\def\mz{\mathbb{Z}}

\def\mf{\mathbb{F}}

\def\Lrho{L_{\widetilde{\rho}}}
\def\pf{\textbf{Proof. }}
\newtheorem{lemma}{Lemma}[section]
\newtheorem{theorem}[lemma]{Theorem}

\begin{document}

\title{Modular Elliptic Curves over the Field of Twelfth Roots of Unity}
\author{Andrew Jones}
\maketitle

\begin{abstract} In this article we perform an extensive study of the spaces of automorphic forms for $\mathrm{GL}_2$ of weight two and level $\mathfrak{n}$, for $\mathfrak{n}$ an ideal in the ring of integers of the quartic CM field $\mq(\zeta_{12})$ of twelfth roots of unity. This study is conducted through the computation of the Hecke module $H^*(\Gamma_0(\mathfrak{n}),\mc)$, and the corresponding Hecke action. Combining this Hecke data with the Faltings-Serre method for proving equivalence of Galois representations, we are able to provide the first known examples of modular elliptic curves over this field.
\end{abstract}

\specialsection{\textbf{Introduction}}

Following the proof of the Shimura-Taniyama-Weil conjecture at the turn of the $21^{st}$ century, it is known that all rational elliptic curves are \it{modular}, in the sense that for almost all primes $\ell$, the representation arising from the action of the absolute Galois group $G_\mq := \mathrm{Gal}(\overline{\mq}/\mq)$ on the $\ell$-adic Tate module of an elliptic curve $E$ defined over $\mq$ is isomorphic to the $\ell$-adic representation attached to a cuspidal modular form $f$ of weight $2$, whose level matches the conductor of $E$.

This notion of modularity can be viewed within the more general framework of the Langlands program, which conjectures in particular that certain $n$-dimensional $\ell$-adic representations of the absolute Galois group $G_F$ (where $F$ is a number field) should correspond to automorphic forms for the algebraic group $\mathrm{Res}_{F/\mq}(\mathrm{GL}_n)$, where $\mathrm{Res}_{F/\mq}$ denotes the Weil restriction of scalars from $F$ to $\mq$. Amongst these Galois representations should be the $2$-dimensional representation arising from the action of $G_F$ on the $\ell$-adic Tate module of an elliptic curve $E$ defined over $F$, so it is natural to consider modularity over an arbitrary number field.

One can approach this question from either a theoretical or a computational viewpoint. On the theoretical side, there is the search for a concrete description of the Galois representation attached to an automorphic form for the group $\mathrm{Res}_{F/\mq}(\mathrm{GL}_2)$ for a general number field $F$. Carayol (\textbf{[Car86]}), Taylor (\textbf{[Tay89]}) and Blasius-Rogawski (\textbf{[BR89]}) defined these representations for the case of Hilbert modular forms over totally real fields whose weights are at least $2$ (including, in particular, those forms which should correspond to elliptic curves). Subsequently, Jarvis and Manoharmayum (\textbf{[JM08]}) proved modularity of semistable elliptic curves over the real quadratic fields $\mq(\sqrt{2})$ and $\mq(\sqrt{17})$, and recent work by Freitas, Le Hung and Siksek (\textbf{[FHS13]}) establishes modularity of elliptic curves over \it{all} real quadratic fields.

Harris, Soudry and Taylor (\textbf{[HST93]} and \textbf{[Tay94]}) provided a construction in the case of automorphic forms over imaginary quadratic fields, which was later refined by Berger and Harcos (\textbf{[BH07]}). More recently, Mok (\textbf{[Mok14]}) has constructed Galois representations attached to automorphic forms defined over CM fields, subject to certain conditions on the central character of these forms (in particular, his construction covers those forms which are expected to correspond to elliptic curves). Harris, Lan, Taylor and Thorne (\textbf{[HLTT13]}) and Scholze (\textbf{[Sch13]}) have removed these restrictions, and in fact construct Galois representations attached to regular automorphic forms for $\mathrm{Res}_{F/\mq}(\mathrm{GL}_n)$ for all $n$ over CM or totally real fields (with Scholze's results extending even further to account for ``torsion automorphic forms''). 

As yet, there are no general modularity results over these latter fields, although examples of modular elliptic curves have been found over certain imaginary quadratic fields, such as in \textbf{[DGP10]}. In this paper, the authors present an algorithm to test for isomorphism of two representations of the absolute Galois group $G_F$, which relies on knowledge of the traces of these representations on Frobenius elements over primes in $F$.

As in the classical case, these traces are determined by the action of Hecke operators on spaces of newforms, which has led to the development of methods for constructing spaces of automorphic forms defined over number fields, and for computing the Hecke action on these spaces. Classically, one can study modular forms through the use of \it{modular symbols} (see, for example, \textbf{[Ste07]}). By the Eichler-Shimura isomorphism, one can realise cuspidal Hecke eigenforms as classes in the cohomology of modular curves, and modular symbols provide us with a means to compute this cohomology, and the corresponding Hecke action. More generally, it is known that cuspidal automorphic forms for $\mathrm{Res}_{F/\mq}(\mathrm{GL}_2)$ can be realised as classes in the cohomology of certain locally symmetric spaces, and so one might hope to extend the notion of modular symbols to a wider variety of number fields. 

Until recently, attention had largely been focused on imaginary quadratic or totally real fields; however, the authors of \textbf{[GHY13]} provide a description of a method for computing automorphic forms over the \it{quartic} CM field $\mq(\zeta_5)$, where $\zeta_5$ denotes a primitive fifth root of unity. Moreover, they are able to compute the action of Hecke operators on these forms, and present several examples of elliptic curves which appear to be modular. In this paper, we shall extend these methods to the field $\mq(\zeta_{12})$, and, adapting the methodology of \textbf{[DGP10]}, shall combine the resulting data with our knowledge of the Galois representations constructed in \textbf{[Mok14]} to give the first proven examples of modular elliptic curves over a quartic CM field.

\newpage\specialsection{\textbf{Cuspidal Automorphic Forms and Representations}}

We begin with a brief discussion of automorphic forms for the group $\mathrm{Res}_{F/\mq}(\mathrm{GL}_2)$, where $F$ is a quartic CM field, before describing the attached Galois representations constructed in \textbf{[Mok14]}.

Let $F$ be such a field, with ring of integers $\mathcal{O}$, let $G = \mathrm{Res}_{F/\mq}(\mathrm{GL}_2)$, and let $G(\mathbb{A}) = G(\mathbb{A}_f) \times G(\mr)$ be a decomposition of the adelic points of $G$ into finite and infinite parts. For simplicity, we shall assume that $F$ has class number $1$. We can identify $G(\mr)$ with the product $\mathrm{GL}_2(\mc) \times \mathrm{GL}_2(\mc)$ via a fixed choice of non-conjugate embeddings of $F$ into the complex numbers, and subsequently define a compact subgroup $K_\infty$ of $G(\mr)$ by $K_\infty = U(2) \times U(2)$. 

For an ideal $\mathfrak{n}$ of $F$, let $K_0(\mathfrak{n})$ denote the product of the local subgroups $K_v(\mathfrak{n})$, where $v$ runs through the non-archimedean places of $F$, and $K_v(\mathfrak{n})$ is equal to either the subgroup of $\mathrm{GL}_2(\mathcal{O}_v)$ comprising those matrices which are upper triangular modulo $\mathfrak{n}$, if $v|\mathfrak{n}$, or else is equal to $\mathrm{GL}_2(\mathcal{O}_v)$. Note that $K_0(\mathfrak{n})$ is a compact open subgroup of $G(\mathbb{A}_f)$.

By an \it{automorphic form of level} $\mathfrak{n}$, we shall mean a function $\varphi: G(\mq)\backslash G(\mathbb{A}) \rightarrow \mc$ which acts smoothly on the real Lie group $G(\mr)$, is invariant under the right regular action of $K_0(\mathfrak{n})$ on $G(\mathbb{A}_f)$, whose translates under the group $K_\infty$ span a finite-dimensional vector space, and which is $\mathcal{Z}$-\it{finite} and \it{of moderate growth}, in the sense of \textbf{[Bum94], Section 3.2}. We say that $\varphi$ is \it{cuspidal} if \begin{equation*}\underset{F\backslash \mathbb{A}_F}{\displaystyle\int} \varphi\left(\left(\begin{smallmatrix}1&x\\0&1\end{smallmatrix}\right)g\right)dx = 0\end{equation*} for all $g \in G(\mathbb{A})$, and has \it{trivial central character}, if it is invariant under the regular action of the centre of $G(\mathbb{A})$.

Such forms lie in the space $L^2_0(G(\mq)\backslash G(\mathbb{A}),\mathbf{1})$ of square integrable functions on $G(\mathbb{A})$ which are invariant under the action of $G(\mq)$, and are cuspidal with trivial central character. Under the right regular action of $G(\mathbb{A})$, this space decomposes into a direct sum of irreducible Hilbert space representations, \begin{equation*}L^2_0(G(\mq)\backslash G(\mathbb{A}),\mathbf{1}) = \underset{\pi}{\bigoplus}~n_\pi V_\pi,\end{equation*} where $n_\pi \in \mathbb{N}$ and $V_\pi$ is a Hilbert space on which $G(\mathbb{A})$ acts by the homomorphism $\pi$, and we define a \it{cuspidal automorphic representation} for $G$ to be any subrepresentation which is isomorphic to one of these summands. Cuspidal automorphic forms of level $\mathfrak{n}$ can be realized in the space of $K_0(\mathfrak{n})$-fixed vectors of some cuspidal automorphic representation, and so we may transfer our attention to these latter objects, which we are better equipped to work with.

It is known that each cuspidal automorphic representation $\pi$ admits a decomposition $\pi = {\otimes}~\pi_v$, where $v$ runs over the places of $F$. For each non-archimedean place $v$, $\pi_v$ is an irreducible \it{admissible} complex representation $\pi_v: \mathrm{GL}_2(F_v) \rightarrow \mathrm{GL}(V_v)$, in the sense that the stabilizer of each point $v \in V_v$ is open, and, given a compact open subgroup $K_v$ of $\mathrm{GL}_2(F_v)$, the space of $K_v$-fixed vectors in $V_v$ is finite-dimensional.

We recall that a representation $\pi_v$ is \it{unramified} if the space of $\mathrm{GL}_2(\mathcal{O}_v)$-fixed vectors is non-trivial. Since $\pi$ is automorphic, each unramified component $\pi_v$ is known to be a \it{principal series}, and is parametrised by a pair $(\chi_1,\chi_2)$ of characters of $F_v^\times$. To each such representation we can assign a semisimple conjugacy class $t_{\pi_v}$ in $\mathrm{GL}_2(\mc)$, known as the \it{Langlands class} of $\pi_v$, a representative of which is given by the matrix $\left(\begin{smallmatrix}\chi_1(\varpi)&0\\0&\chi_2(\varpi)\end{smallmatrix}\right)$, where $\varpi$ is a uniformiser for $\mathcal{O}_v$.

As intimated previously, we can realise automorphic forms (or, as we shall consider, the corresponding automorphic representations) as classes in the cohomology of certain locally symmetric spaces. An in-depth discussion of this can be found in \textbf{[Schw06], Chapter 3}, of which we shall present a brief pr\'{e}cis, tailored to our particular case.

For a fixed ideal $\mathfrak{n}$ in $F$, we define the locally symmetric space \begin{equation*}X_0(\mathfrak{n}) = A_G^0(\mr)G(\mq)\backslash G(\mathbb{A})/K_\infty K_0(\mathfrak{n}),\end{equation*} where the \it{split component} $A_G^0(\mr)$ is isomorphic to $\mr_+$, embedded diagonally into the two components of $G(\mr) \simeq \mathrm{GL}_2(\mc) \times \mathrm{GL}_2(\mc)$. We then define $H^*(X_0(\mathfrak{n}),\mc)$ to be the cohomology of $\Omega(X_0(\mathfrak{n}),\mc)$, where the latter denotes the de Rham complex of smooth, complex-valued differentials on $X_0(\mathfrak{n})$. 

Let $\mathfrak{m}_G$ be the Lie algebra of $A_G^0(\mr)\backslash G(\mr)$, and denote by $\mathcal{A}(K_0(\mathfrak{n}))$ the space of automorphic forms of level $\mathfrak{n}$ on which $A_G^0(\mr)$ acts trivially. There is then an isomorphism of $K_0(\mathfrak{n})$-modules: \begin{equation*}H^*(X_0(\mathfrak{n}),\mc) \simeq H^*(\mathfrak{m}_G,K_\infty;\mathcal{A}(K_0(\mathfrak{n}))),\end{equation*} where $H^*(\mathfrak{m}_G,K_\infty;\mathcal{A}(K_0(\mathfrak{n})))$ denotes the Lie algebra cohomology with respect to $(\mathfrak{m}_G,K_\infty)$ (in the sense of \textbf{[Schw06], Section 3.2}).

There is a decomposition \begin{equation*}H^*(X_0(\mathfrak{n}),\mc) = H^*_{\mathrm{Eis}}(X_0(\mathfrak{n}),\mc) \oplus H^*_{\mathrm{cusp}}(X_0(\mathfrak{n}),\mc)\end{equation*} of $H^*(K_0(\mathfrak{n}),\mc)$ into \it{Eisenstein} and \it{cuspidal} parts, the latter of which is connected to cuspidal automorphic representations. Denoting by $V_\pi = V_{\pi_f} \otimes V_{\pi_\infty}$ the space on which a cuspidal automorphic representation $\pi$ acts, we have a decomposition of $K_0(\mathfrak{n})$-modules \begin{equation*}H^*_{\mathrm{cusp}}(X_0(\mathfrak{n}),\mc) \simeq \underset{\pi}{\displaystyle\bigoplus} ~H^*(\mathfrak{m}_G,K_\infty;V_{\pi_\infty}) \otimes V_{\pi_f}^{K_0(\mathfrak{n})},\end{equation*} where the sum ranges over all cuspidal automorphic subrepresentations of the space $\mathcal{A}(K_0(\mathfrak{n}))$ with trivial central character, and $V_{\pi_f}^{K_0(\mathfrak{n})}$ denotes the space of $K_0(\mathfrak{n})$-fixed vectors in $V_{\pi_f}$ (see \textbf{[Schw06], Theorem 4.1}). We shall say that a representation $\pi$ is \it{of cohomological type} and \it{weight two} if $H^*(\mathfrak{m}_G,K_\infty;V_{\pi_\infty})$ is non-zero.

\newpage With this in mind, we can now state a version of the main result in \textbf{[Mok14]}:

\begin{theorem} Let $F$ be a CM field, and let $\pi$ be a cuspidal automorphic representation of $\mathrm{Res}_{F/\mq}(\mathrm{GL}_2)$ of cohomological type, with trivial central character, and fix a prime $\ell$. Then there exists an $\ell$-adic Galois representation \begin{equation*}\rho_\pi: G_F \rightarrow \mathrm{GL}_2(\overline{\mq}_\ell)\end{equation*} such that, for each place $v$ of $F$ not dividing $\ell$, we have the local-to-global compatibility statement, up to semisimplification:
\begin{equation*}\mathrm{WD}(\rho_{\pi,v})^{ss} \simeq \mathscr{L}_v(\pi_v \otimes |\mathrm{det}|_v^{-\frac{1}{2}})^{ss}.\end{equation*}
Furthermore, if $\pi_v$ is not a twist of Steinberg (e.g., is an unramified principal series) then we have the full local-to-global compatibility statement, up to Frobenius semisimplification:
\begin{equation*}\mathrm{WD}(\rho_{\pi,v})^{\mathrm{Frob}} \simeq \mathscr{L}_v(\pi_v \otimes |\mathrm{det}|_v^{-\frac{1}{2}}).\end{equation*}
\label{Mokrep}\end{theorem}

Here $\mathcal{L}_v(\pi_v)$ denotes the representation of the Weil-Deligne group $W'_v$ of $F_v$ assigned to $\pi_v$ by the local Langlands correspondence for $\mathrm{GL}_2$, and $\mathrm{WD}$ denotes the Weil-Deligne functor taking representations of $G_{F_v}$ to representations of $W'_v$.

We shall say that an elliptic curve $E$ over $F$ is \it{modular} if, for some rational prime $\ell$, the $\ell$-adic Galois representation $\rho_E$ defined by the action on the $\ell$-adic Tate module of $E$ is, up to semisimplification, isomorphic to $\rho_\pi$ for some cohomological cuspidal automorphic representation $\pi$ of weight two and level $\mathfrak{n}_E$, where $\mathfrak{n}_E$ denotes the conductor of $E$. Note that, since such representations form a compatible family ranging over all rational primes $\ell$, isomorphism for some prime implies isomorphism for all but finitely many primes.

We note some important details concerning these representations. Firstly, since the Frobenius semisimplification of a Weil-Deligne representation agrees with the original representation on inertia, $\mathrm{WD}(\rho_{\pi,v})$ is unramified at a place $v$ of $F$ if, and only if, $\rho_{\pi,v}$ is too. Similarly, under the local Langlands correspondence, $\mathcal{L}_v(\pi_v)$ and $\pi_v$ are unramified at the same set of places. In particular, if the representation $\pi$ occurs in the cohomology of the locally symmetric space $X_0(\mathfrak{n})$, then $\rho_{\pi,v}$ is unramified at all primes not dividing $\mathfrak{n}$.

Secondly, the determinant of $\mathcal{L}_v(\pi_v)$ corresponds to the central character of $\pi_v$, and thus is trivial, so in particular the determinant of $\rho_{\pi,v}$ is equivalent to $|\mathrm{det}|_v^{-\frac{1}{2}}$. Observing that, for a uniformiser $\varpi$, we have $\mathrm{det}(\rho_{\pi,v}(\mathrm{Frob_v})) = |\mathrm{det}(\varpi)|_v^{-\frac{1}{2}} = q$, where $q$ is the cardinality of the residue field of $F_v$, we see that $\mathrm{det}(\rho_{\pi,v})$ is equivalent to the local cyclotomic character.

Finally, the trace of $\mathcal{L}_v(\pi_v)(\mathrm{Frob_v})$ is equal to the trace of the Langlands class $t_{\pi_v}$ of $\pi_v$, and so $\mathrm{Tr}(\rho_{\pi,v}(\mathrm{Frob_v})) = q^\frac{1}{2}\mathrm{Tr}(t_{\pi_v})$. If the representation $\pi_v$ is unramified, the space of $\mathrm{GL}_2(\mathcal{O}_v)$-fixed vectors of $\pi_v$ is known to be one-dimensional, on which the Hecke operator $T_v$ (defined as the  normalized characteristic function of the double coset $\mathrm{GL}_2(\mathcal{O}_v)\left(\begin{smallmatrix}\varpi&0\\0&1\end{smallmatrix}\right)\mathrm{GL}_2(\mathcal{O}_v)$) acts via scalar multiplication by precisely this value.

\specialsection{\textbf{Comparing Galois Representations}}

Armed with our knowledge of these representations, we can use the following method, first exhibited in \textbf{[DGP10]} for imaginary quadratic fields, to determine whether an elliptic curve $E$ defined over a quartic CM field is modular. Throughout, we shall denote by $\rho_E$ and $\rho_\pi$ the Galois representations corresponding to the curve $E$ and a cuspidal automorphic representation $\pi$ respectively.

To begin with, we need to consider the \it{residual representation} $\widetilde{\rho}$ of an arbitrary $\ell$-adic Galois representation $\rho$, whose definition we now recall. It is a standard result that, up to isomorphism, the image of $\rho$ can be assumed to be defined over some finite extension of $\mz_\ell$, and thus we may compose $\rho$ with a reduction map to the residue field of this valuation ring to obtain a representation $\bar{\rho}: G_F \rightarrow \mathrm{GL}_2(\mf_{\ell^r})$ for some $r$. We define $\widetilde{\rho}: G_F \rightarrow \mathrm{GL}_2(\mf_{\ell^r})$ to be the semisimplification of any such $\bar{\rho}$ (which is well-defined up to isomorphism). 
 
Each Galois representation that we shall consider is \it{rational}, in the sense that, for each non-archimidean place $v$, the coefficients of the characteristic polynomial of $\mathrm{Frob}_v$ are rational. The following result (which is no doubt well known to the experts, but for which we could find no reference) shows that, for any such representation, the corresponding residual representation can be assumed to have image in $\mathrm{GL}_2(\mf_\ell)$.

\begin{lemma} Let $\rho: G_F \rightarrow \mathrm{GL}_2(\overline{\mq}_\ell)$ be a rational Galois representation, and let $\widetilde{\rho}: G_F \rightarrow \mathrm{GL}_2(\overline{\mf}_\ell)$ be the corresponding residual representation. Then there exists an element $t \in \mathrm{GL}_2(\overline{\mf}_\ell)$ such that $t\widetilde{\rho}(g)t^{-1} \in \mathrm{GL}_2(\mf_\ell)$ for all $g \in G_F$.\end{lemma}

\pf Let the image of $\widetilde{\rho}$ lie in $\mathrm{GL}_2(\mf_{\ell^{r}})$, for some $r$, let $\sigma$ be a generator of the cyclic group $\mathrm{Gal}(\mf_{\ell^{r}}/\mf_\ell)$, and consider the representations $\widetilde{\rho}$ and $\sigma \circ \widetilde{\rho}$. Since $\sigma$ fixes $\mf_\ell$, the characteristic polynomials of these two representations on Frobenius elements coincide, and, since both representations are semisimple by definition, a theorem of Brauer and Nesbitt implies that they are isomorphic. Thus there exists some $s \in \mathrm{GL}_2(\mf_{\ell^{r}})$ such that $\sigma(\widetilde{\rho}(g)) = s^{-1}\widetilde{\rho}(g)s$ for all $g \in G_F$.

Define an element \begin{equation*}\pi_s = \mathrm{Nm}_{\mf_{\ell^{r}}/\mf_\ell}(s) = \displaystyle\prod_{i=1}^{r}\sigma^i(s) \in \mathrm{GL}_2(\mf_\ell).\end{equation*} Since $\mathrm{GL}_2(\mf_\ell)$ is finite, $\pi_s$ must have finite order, $m$, say. Letting $n = mr$, we find that $\mathrm{Nm}_{\mf_{\ell^{n}}/\mf_\ell}(s) = \mathrm{Nm}_{\mf_{\ell^{n}}/\mf_{\ell^{r}}}(\pi_s) = \pi_s^m = \mathrm{Id}$.

Now, let $\tau$ be a generator of the cyclic group $G = \mathrm{Gal}(\mf_{\ell^{n}}/\mf_\ell)$, so that in particular $\tau|_{\mf_{\ell^{r}}} = \sigma$. One can define a $G$-cocycle $\gamma$ by the map \begin{equation*}\gamma: G \rightarrow \mathrm{GL}_2(\mf_{\ell^{n}});~\tau^k \mapsto \gamma_{\tau^k} := \underset{i=0}{\overset{k-1}{\displaystyle\prod}} ~\tau^i(s).\end{equation*} By Hilbert's Theorem 90, the first cohomology group $H^1(G,\mathrm{GL}_2(\mf_{\ell^{n}}))$ is trivial, and thus there exists some $t \in \mathrm{GL}_2(\mf_{\ell^{n}})$ such that $t\gamma_{\tau^k}\tau(t)^{-k} = \mathrm{Id}$ for all $k$. In particular, since $\gamma_\tau = s$, we have $s = t^{-1}\tau(t)$.

Viewing $\mathrm{GL}_2(\mf_{\ell^{r}})$ as a subgroup of $\mathrm{GL}_2(\mf_{\ell^{n}})$ on which $\tau \in \mathrm{Gal}(\mf_{\ell^{n}}/\mf_{\ell^{r}})$ acts as the element $\sigma \in \mathrm{Gal}(\mf_{\ell^{r}}/\mf_\ell)$, we observe that, for each $g \in G_F$, \begin{equation*}\tau(t\widetilde{\rho}(g)t^{-1}) = \tau(t)\sigma(\widetilde{\rho}(g))\tau(t)^{-1} = \tau(t)s^{-1}\widetilde{\rho}(g)s\tau(t)^{-1} = t\rho(g)t^{-1},\end{equation*} and thus $t\rho(g)t^{-1} \in \mathrm{GL}_2(\mf_\ell)$ for all $g \in G_F$, as required. \hfill $\square$

\textbf{Remark:} There is nothing special about our restriction to rational representations. The same argument shows that, if the coefficients of the characteristic polynomial of $\widetilde{\rho}(\mathrm{Frob}_\mathfrak{p})$ lie in $\mf_{\ell^r}$ for all primes $\mathfrak{p}$, then we can define $\widetilde{\rho}$ over the field $\mf_{\ell^r}$.

We restrict our attention to the case $\ell = 2$. Thus, for any rational $2$-adic representation $\rho$, the residual representation $\widetilde{\rho}$ has image in $\mathrm{GL}_2(\mf_2)$, which is isomorphic to $S_3$. It is straightforward to see that, up to semisimplification, $\widetilde{\rho}$ has only three possible images; it can be trivial, have cyclic image isomorphic to $C_3$, or be isomorphic to $S_3$ itself.

Using class field theory, it is possible to determine isomorphism of the residual representations $\widetilde{\rho_E}$ and $\widetilde{\rho_\pi}$ by checking certain parity conditions on the traces of these representations at Frobenius elements above a finite set $\mathcal{S}_1$ of prime ideals of $F$ (we discuss this method in \textbf{Appendix A}). Henceforth, we shall assume that these residual representations are indeed isomorphic. It is then known that one can produce a second set $\mathcal{S}_2$ of prime ideals of $F$ such that the full representations $\rho_E$ and $\rho_\pi$ are isomorphic if, and only if, the traces of $\rho_E(\mathrm{Frob}_\mathfrak{p})$ and $\rho_\pi(\mathrm{Frob}_\mathfrak{p})$ for all primes $\mathfrak{p} \in \mathcal{S}_2$ are equal. A suitable set $\mathcal{S}_2$ can be constructed by the \it{Faltings-Serre method} in the case that the residual representations are absolutely irreducible, or by a result of Livn\'{e} in all other cases.

For reference, we state these results here, beginning with the absolutely irreducible case:

\begin{theorem} Let $\rho_1$, $\rho_2: G_F \rightarrow \mathrm{GL}_2(\mz_2)$ be two representations which have the same determinant, are unramified outside a finite set of primes in $F$, and whose residual representations are absolutely irreducible and isomorphic. Then there exists a finite set $\mathcal{S}$ of primes in $F$ such that $\rho_1$ and $\rho_2$ have isomorphic semisimplifications if, and only if, $\mathrm{Tr}(\rho_1(\mathrm{Frob}_\mathfrak{p})) = \mathrm{Tr}(\rho_2(\mathrm{Frob}_\mathfrak{p}))$ for all $\mathfrak{p} \in \mathcal{S}$.\end{theorem}

For a full, constructive account of the proof, we refer the reader to \textbf{[DGP10], Section 4}. The requirement that both representations have image in $\mathrm{GL}_2(\mz_2)$, which \it{a priori} could be problematic, is in fact easily dealt with. Indeed, absolute irreducibility of the residual representations, combined with the fact that the traces of the full representations are rational, means that we can use a result of Carayol (see \textbf{[Car94], Theorem 2}) to find that $\rho_\pi$ is in fact equivalent to a representation which takes values in $\mathrm{GL}_2(\mz_2)$ (of course, we already know this to be the case for $\rho_E$).

If the images of $\rho_\pi$ and $\rho_E$ are not absolutely irreducible, then we require the following result of Livn\'{e}, whose statement we borrow from the thesis of Ch\^{e}nevert (\textbf{[Ch\^{e}08], Theorem 5.4.9}):

\begin{theorem}Let $F$ be a number field, and $V_\lambda$ a finite extension of $\mq_2$ with ring of integers $\mathcal{O}_\lambda$ and maximal ideal $\lambda$. Let \begin{equation*}\rho_1,\rho_2:G_F \rightarrow \mathrm{GL}_2(V_\lambda)\end{equation*} be two continuous representations unramified outside a finite set $S$ of places of $F$, such that \begin{equation*}\mathrm{Tr}(\widetilde{\rho_1}) \equiv \mathrm{Tr}(\widetilde{\rho_2}) \equiv 0 ~(\mathrm{mod}~\lambda),~\mathrm{and}~\mathrm{det}(\widetilde{\rho_1}) \equiv \mathrm{det}(\widetilde{\rho_2}) \equiv 1 ~(\mathrm{mod}~\lambda).\end{equation*}

Let $F_{2,S}$ denote the compositum of all quadratic extensions of $F$ unramified outside $S$, and suppose there exists a set of prime ideals $T$ of $\mathcal{O}_F$, disjoint from $S$, such that:
\begin{enumerate}[\upshape(i)]
\item $\{\mathrm{Frob}_\mathfrak{p},\mathfrak{p} \in T\}$ surjects onto $\mathrm{Gal}(F_{2,S}/F)$; and
\item The characteristic polynomials of $\rho_1$ and $\rho_2$ at the elements $\{\mathrm{Frob}_\mathfrak{p},\mathfrak{p} \in T\}$ are equal.
\end{enumerate}

Then $\rho_1$ and $\rho_2$ have isomorphic semisimplifications. \label{Livthm}\end{theorem}

If the residual representations are trivial, then this suffices to prove isomorphism (up to semisimplification) of the full representations. However, if the residual representations have $C_3$-image, then the hypotheses of the theorem do not hold (since the traces of the elements in $\mathrm{GL}_2(\mf_2)$ of order $3$ are odd), and we need an additional step.

Let $\widetilde{\rho_E} \simeq \widetilde{\rho_\pi} \simeq \widetilde{\rho}$, say, and define the fixed field $L = \overline{F}^{\mathrm{ker}(\widetilde{\rho})}$, so that $L/F$ is a cubic Galois extension, and $G_L$ is a normal subgroup of $G_F$. Let $\rho'_E$ and $\rho'_\pi$ denote the restrictions of our original representations to $G_L$. Then it is clear that the corresponding residual representations are trivial, and we can use \textbf{Theorem \ref{Livthm}} to determine isomorphism of $\rho'_E$ and $\rho'_\pi$ (as usual, up to semisimplification). 

For the following argument, we impose an additional restriction on the elliptic curve $E$, namely that the base change of $E$ to the field $L$ does not possess complex multiplication. If this restriction is satisfied, then the representation $\rho'_E$ (and thus also $\rho'_\pi$) is in fact irreducible. By Schur's lemma, $\mathrm{Hom}_{G_L}(\rho'_E,\rho'_\pi)$ contains a copy of $V_\lambda$ on which $G_L$ acts trivially, and so $\mathrm{Hom}_{G_L}(\mathbf{1},(\rho_E \otimes \rho_\pi^\vee)|_{G_L})$ is non-trivial. By Frobenius reciprocity, the latter group is isomorphic to $\mathrm{Hom}_{G_F}(\mathrm{Ind}_{G_L}^{G_F}(\mathbf{1}), \rho_E \otimes \rho_\pi^\vee)$, which decomposes as a direct sum \begin{equation*}\mathrm{Hom}_{G_F}(\mathrm{Ind}_{G_L}^{G_F}(\mathbf{1}), \rho_E \otimes \rho_\pi^\vee) \simeq \underset{\chi|_{G_L} = \mathbf{1}}{\bigoplus} \mathrm{Hom}_{G_F}(\rho_\pi \otimes \chi, \rho_E).\end{equation*}
Invoking Schur's lemma once more, we observe that one of these summands must be non-trivial, and that $\rho_E \simeq \rho_\pi \otimes \chi$ for some character $\chi$ of $G_F$ whose restriction to $G_L$ is non-trivial. One can then determine whether this character is trivial, by finding a prime $\mathfrak{p}$ of $F$ which is inert in $L$. In this case, $\mathrm{Frob}_\mathfrak{p}$ is non-trivial, and so $\chi$ is completely determined by the value it takes on this Frobenius element. In particular, if $\mathrm{Tr}(\rho_\pi(\mathrm{Frob}_\mathfrak{p})) = \mathrm{Tr}(\rho_E(\mathrm{Frob}_\mathfrak{p}))$, then $\chi(\mathrm{Frob}_\mathfrak{p}) = 1$, $\chi$ is trivial, and $\rho_\pi$ and $\rho_E$ have isomorphic semisimplifications, as required.

\specialsection{\textbf{Computing the Hecke Action on Cohomology}}

To utilise the method described in the previous section, we require knowledge of the representation $\rho_E$, and the values $\mathrm{Tr}(\rho_\pi(\mathrm{Frob}_\mathfrak{p}))$ for a finite set of primes $\mathfrak{p}$ of $F$. The former is straightforward to find, given an elliptic curve $E$, as the action of $G_F$ on the $2$-adic Tate module $\mathrm{Ta}_2(E)$ of $E$ is well understood. We therefore require a method to compute the trace of $\rho_\pi$ on Frobenius elements.

From \textbf{Section 2}, we know that, for a place $v$ of $F$, dividing $\mathfrak{p}$, at which $\pi$ is unramified, this is given by the action of the Hecke operator $T_v$ on the space of $\mathrm{GL}_2(\mathcal{O}_v)$-fixed vectors of $\pi_v$. Moreover, this Hecke action translates to the cohomological setting mentioned previously, where it is described in terms of the action of a double coset operator on the de Rham cohomology of differential forms on $X_0(\mathfrak{n})$.

At this juncture, we make note of an important observation: one could, instead, work with the \it{rational} cohomology, $H^*_{\mathrm{cusp}}(X_0(\mathfrak{n}),\mq)$. By \textbf{[Har06], Chapter III, Proposition 2.2}, this is also a $K_0(\mathfrak{n})$-module, and so admits a Hecke action. Moreover, the identification \begin{equation*}H^*_{\mathrm{cusp}}(X_0(\mathfrak{n}),\mc) \simeq H^*_{\mathrm{cusp}}(X_0(\mathfrak{n}),\mq) \otimes \mc,\end{equation*} is an isomorphism of Hecke modules, and so in particular the action of the Hecke operators $T_v$ on the \it{complex} cohomology can be defined rationally, a fact which shall prove useful later.

We shall now provide a more concrete realization of this cohomology and the corresponding Hecke operators, which we can use to compute the desired information about the Galois representations $\rho_\pi$. Key to this is the fact that we may reinterpret the locally symmetric space $X_0(\mathfrak{n})$ as a quotient of some \it{globally} symmetric space $X$ by an arithmetic subgroup $\Gamma_0(\mathfrak{n})$ of $G(\mq)$, which in turn allows us to compute $H^*(X_0(\mathfrak{n}),\mc)$ by considering the group cohomology $H^*(\Gamma_0(\mathfrak{n}),\mc)$.

More precisely, define the globally symmetric space \begin{equation*}X = G(\mr)/A_G^0(\mr)K_\infty,\end{equation*} and let $\Gamma_0(\mathfrak{n})$ denote the arithmetic subgroup of $G(\mq)$ - which we identify with $\mathrm{GL}_2(F)$ - comprising those matrices in $\mathrm{GL}_2(\mathcal{O})$ which are upper-triangular modulo the ideal $\mathfrak{n}$. Recalling that we have restricted our attention to fields $F$ with trivial class group, we have an identification \begin{equation*}X_0(\mathfrak{n}) \simeq \Gamma_0(\mathfrak{n})\backslash X.\end{equation*} 

A standard argument (an example of which can be found in the appendix of \textbf{[Hid93]}) then states that we have an equivalence of cohomology groups: \begin{equation*}H^*(X_0(\mathfrak{n}),\mc) \simeq H^*(\Gamma_0(\mathfrak{n}),\mc).\end{equation*}
 
We may therefore choose to work with the group cohomology appearing on the right-hand side, but we require a little extra work before we can compute it effectively. 

\newpage Recall that the \it{cohomological dimension} $\nu$ of a torsion-free arithmetic subgroup $\Gamma$ of $G(\mq)$ is defined to be the smallest integer such that $H^{\nu+1}(\Gamma,\mathcal{M}) = 0$ for \it{all} coefficient systems $\mathcal{M}$. Since $\Gamma_0(\mathfrak{n})$ is not torsion-free, we need a more general definition, that of the \it{virtual cohomological dimension} (which we also denote by $\nu$) of an arbitrary arithemtic subgroup $\Gamma$ of $G(\mq)$, which is defined to be the cohomological dimension of any finite index torsion-free subgroup of $\Gamma$. 

\it{Borel-Serre duality} (see \textbf{[BS73], Section 11.4}) then states that we have an equivalence of Hecke modules \begin{equation*}H^{\nu-k}(\Gamma_0(\mathfrak{n}),\mc) \simeq H_k(\Gamma_0(\mathfrak{n}),\mathrm{St}_2 \otimes_\mz \mc),\end{equation*} where the \it{Steinberg module} $\mathrm{St}_2$ for $\mathrm{GL}_2(F)$ is defined by the short exact sequence \begin{equation*}0 \rightarrow \mathrm{St}_2 \rightarrow \mz[\mathbb{P}^1(F)] \stackrel{\epsilon}{\rightarrow} \mz \rightarrow 0\end{equation*} (here $\epsilon$ denotes the augmentation map sending $\sum n_P P$ to $\sum n_P$).

The group homology of $\Gamma_0(\mathfrak{n})$ with coefficients in $\mathrm{St}_2 \otimes_\mz \mc$ can in turn be computed by constructing an appropriate resolution of the Steinberg module. Our preferred example of such a resolution is given by the \it{sharbly complex} $\mathscr{S}_*$, which we define as follows: for each non-negative integer $k$, let $\mathscr{S}_k$ denote the space of $\mz$-linear combinations of $k$-\it{sharblies} -- $(k+2)$-tuples $\mathbf{u} = [u_1,\ldots,u_{k+2}]$ with each $u_i \in \mathcal{O}^2$ -- subject to the relations:
\begin{itemize}
\item $[u_1,\ldots,u_{k+2}] = \mathrm{sgn}(\sigma)[u_{\sigma(1)},\ldots,u_{\sigma(k+2)}]$ for any permutation $\sigma \in S_{k+2}$;
\item $[u,u_2,\ldots,u_{k+2}] = [v,u_2,\ldots,v_{k+2}]$ if there exists $\lambda \in \mr_+$ such that for each embedding $\iota: F \hookrightarrow \mc$, we have $\iota(uu^*) = \lambda \iota(vv^*)$; 
\item $[u_1,\ldots,u_{k+2}] = 0$ if $u_1, \ldots, u_{k+2}$ span a $1$-dimensional $F$-vector space (we call such sharblies \it{degenerate}).
\end{itemize}

To give $\mathscr{S}_*$ the form of a simplicial complex, we equip it with the boundary map $\partial_k: \mathscr{S}_k \rightarrow \mathscr{S}_{k-1}$ defined by \begin{equation*}\partial_k [u_1,\ldots,u_{k+2}] = \underset{i=1}{\overset{k}{\displaystyle\sum}}~(-1)^{i+1}[u_1,\ldots,\hat{u}_i,\ldots,u_{k+2}],\end{equation*} where $\hat{u}_i$ indicates that we omit the vector $u_i$ from the resulting sharbly. The sharbly complex admits an obvious action of $\mathrm{GL}_2(\mathcal{O})$, given by \begin{equation*}g\cdot [u_1,\ldots,u_{k+2}] = [gu_1,\ldots,gu_{k+2}], ~g \in \mathrm{GL}_2(\mathcal{O}),\end{equation*} which commutes with the boundary map, and thus in particular we can define the subcomplex of $\Gamma_0(\mathfrak{n})$-invariants, which we denote by $(\mathscr{S}_*)_{\Gamma_0(\mathfrak{n})}$, by imposing the additional relation that \begin{itemize}
\item $[u_1,\ldots,u_{k+2}] = \gamma [u_1,\ldots,u_{k+2}]$ for all $\gamma \in \Gamma_0(\mathfrak{n})$.
\end{itemize}

It is known (see, for example, \textbf{[AGM11], Theorem 5}) that the sharbly complex provides an acyclic resolution of the Steinberg module, and thus we have an identification \begin{equation*}H_k(\Gamma_0(\mathfrak{n}),\mathrm{St}_2 \otimes_\mz \mc) \simeq H_k((\mathscr{S}_*)_{\Gamma_0(\mathfrak{n})},\mc)\end{equation*} of homology groups. Moreover, the identifications we have established are all Hecke-equivariant (see, for example, \textbf{[AGM13], Theorem 2.4}), and so we have a Hecke action on the sharbly homology. One can see without too much difficulty that, for a non-archimedean place $v$ of $F$, corresponding to the prime ideal $\mathfrak{p}$ of $F$, the action of the Hecke operator $T_v$ on $\mathscr{S}_*$ is given by \begin{equation*}T_v(u) = \displaystyle\sum~g_i\mathbf{u},\end{equation*} where we have a decomposition of the double coset space \begin{equation*}\Gamma_0(\mathfrak{n})\left(\begin{smallmatrix}1&0\\0&\nu\end{smallmatrix}\right)\Gamma_0(\mathfrak{n}) = \displaystyle\coprod~\Gamma_0(\mathfrak{n})g_i,\end{equation*} for a generator $\nu$ of the ideal $\mathfrak{p}$.

Thus we can, in theory at least, use the homology of the sharbly complex to compute the Hecke action on cuspidal automorphic forms of level $\mathfrak{n}$. Moreover, as we are interested in classes which correspond to cuspidal automorphic forms, we can restrict our attention to homology whose degree lies within a specified range. 

Indeed, one can show (using a slight adaptation of \textbf{[Schw06], Theorem 6.2}) that the cuspidal cohomology $H^i_{\mathrm{cusp}}(X_0(\mathfrak{n}),\mc)$ is non trivial only if $2 \leq i \leq 5$. Combined with the fact that the virtual cohomological dimension of $\Gamma_0(\mathfrak{n})$ is $6$ (which can be found by applying \textbf{[BS73], Theorem 11.4.4} to our specific case) means that we need only consider homology in degrees $1$ to $4$. Furthermore, it is known that any cuspidal automorphic form which appears as a homology class in one of these degrees in fact appears in \it{every} degree, so we lose nothing by specializing to a single degree. 

The sharbly complex is intrinsically linked with the symmetric space $X$ defined previously, and so we shall spend some time discussing the geometry of this space in greater detail. Recall that \begin{equation*}X = G(\mr)/A_G^0(\mr)K_\infty,\end{equation*} and moreover that $G(\mr) \simeq \mathrm{GL}_2(\mc) \times \mathrm{GL}_2(\mc)$, $K_\infty \simeq \mathrm{U}(2) \times \mathrm{U}(2)$, and $A_G^0(\mr) \simeq \mr_+$. Using the standard identification of $\mathrm{SL}_2(\mc)/\mathrm{SU}(2)$ with the hyperbolic $3$-space $\mathcal{H}_3$, we therefore find that \begin{equation*}X \simeq \mathcal{H}_3 \times \mathcal{H}_3 \times \mr_+,\end{equation*} a $7$-dimensional space.

We can identify $X$ with a cone of binary Hermitian forms. More precisely, let $v_1$ and $v_2$ be two non-conjugate embeddings of $F$ into $\mc$, and define $\mathrm{Herm}_2(F_{v_i})$ to be the space of $2 \times 2$ Hermitian matrices with entries in $F_{v_i}$. Then we can define an inner product space \begin{equation*}\mathcal{V} = \mathrm{Herm}_2(F_{v_1}) \times \mathrm{Herm}_2(F_{v_2}),\end{equation*} with inner product $\langle ~,~\rangle$ defined by \begin{equation*}\langle \Phi,\Psi\rangle = 2\mathrm{Tr}(\Phi_1\Psi_1+\Phi_2\Psi_2)\end{equation*} (the factor of $2$ is largely irrelevant in this case, but comes into play if we extend these ideas to fields of mixed signature -- see \textbf{[GY13]} for an example of this). 

\newpage Note that any point $\Phi \in \mathcal{V}$ defines a Hermitian form on $F^2$ (which, through abuse of notation, we shall also denote by $\Phi$) given by \begin{equation*}\Phi(x,y) = 2\mathrm{Tr}(\Phi_1x_1y_1^* + \Phi_2x_2y_2^*),\end{equation*} where $x_i$ and $y_i$ denote the images of $x$ and $y$ under the embedding $v_i$, and $^*$ denotes the complex conjugate transpose. We shall henceforth refer to points in $\mathcal{V}$ and their corresponding Hermitian forms interchangeably, dependent on the context.

We denote by $\mathcal{C}$ the cone of positive definite forms in $\mathcal{V}$, whose closure $\overline{\mathcal{C}}$ is the cone of positive semi-definite forms. It can be shown that a form $\Phi$ belongs to $\mathcal{C}$ (respectively $\overline{\mathcal{C}}$) if, and only if, each matrix $\Phi_i$ is positive definite (respectively positive semi-definite). The group $G(\mr)$ acts on $\mathcal{V}$ via the map \begin{equation*}g \cdot \Phi = (g_1\Phi_1g_1^*,g_2\Phi_2g_2^*),\end{equation*} and this action preserves $\mathcal{C}$. Moreover, it can easily be seen that $G(\mr)$ acts \it{transitively} on $\mathcal{C}$, and that the stabilizer of any point is isomorphic to $K_\infty$. Thus we obtain an isomorphism \begin{equation*}\mathcal{C}/\mr_+ \simeq X.\end{equation*}

The usefulness of this isomorphism lies in the fact that the cone $\mathcal{C}$ is an example of a \it{positivity domain}, in the sense of Koecher's work in \textbf{[Koe60]}. Such structures possess many desirable properties: in particular, $\mathcal{C}$ admits a decomposition into polyhedral cones, which will provide us with a cellular decomposition of our symmetric space $X$.

This decomposition is straightforward to describe. To each vector $x \in \mathcal{O}^2$, we can define a point $q(x) \in \mathcal{V}$ by setting \begin{equation*}q(x) = (x_1x_1^*, x_2x_2^*).\end{equation*} Note that, using this definition, we can rewrite the second relation defining the sharbly complex as 
\begin{itemize}\item $[u,u_2,\ldots,u_{k+2}] = [v,u_2,\ldots,v_{k+2}]$ if $q(u) = \lambda q(v)$ for some $\lambda \in \mr_+$.\end{itemize} 

Each matrix $x_ix_i^*$ is positive semi-definite, and so in particular $q(x)$ lies in $\overline{\mathcal{C}}$. The set $\Xi := \{q(x);~x \in \mathcal{O}^2\}$ is an example of what Koecher refers to as an \it{admissible subset} of $\overline{\mathcal{C}}$. For each form $\Phi \in \mathcal{C}$, the set of values $\{\langle \Phi, q(x)\rangle:~q(x) \in \Xi\}$ forms a discrete subset of $\mr$, and in particular we may define the \it{minimum} of $\Phi$ to be \begin{equation*}m(\Phi) := \underset{q(x) \in \Xi}{\mathrm{inf}}~\{\langle \Phi,q(x)\rangle\}.\end{equation*} The points $q(x)$ in $\Xi$ for which $\langle \Phi,q(x)\rangle$ attains this minimum are called \it{minimal vectors} of $\Phi$, and we denote by $M(\Phi)$ the set of minimal vectors of $\Phi$. If the set of minimal vectors of a form $\Phi$ span the vector space $\mathcal{V}$, then we call $\Phi$ \it{perfect}. Note that, if $\Phi$ is perfect, then so too is $\lambda \Phi$ for any $\lambda \in \mr_+$, so we may assume without loss of generality that $m(\Phi) = 1$ for all perfect forms $\Phi$.

\newpage The perfect forms provide us with our desired decomposition. To each perfect form $\Phi$, we can assign a convex polytope $\mathcal{F}_\Phi$ in $\overline{\mathcal{C}}$ by taking the convex hull of the minimal vectors of $\Phi$. We can then define a polyhedral cone, known as a \it{perfect pyramid}, by taking the cone above $\mathcal{F}_\Phi$ (that is, the set of half-lines passing through both the origin and a point in $\mathcal{F}_\Phi$). Koecher's work then shows that the set of perfect pyramids provides a decomposition of the cone $\mathcal{C}$, in the sense that every point in $\mathcal{C}$ lies in some perfect pyramid, and any two perfect pyramids have disjoint interiors. Moreover, the action of $\mathrm{GL}_2(\mathcal{O})$ on $\mathcal{O}^2$ clearly induces an action on the set of points $\Xi$, which preserves the set of perfect pyramids, and under this action there are only finitely many equivalence classes of perfect pyramids. For a more detailed description of these ideas, we refer the reader to the excellent exposition in \textbf{[Gun11], Section 7}.

Recalling our earlier identification of the symmetric space $X$ with $\mathcal{C}/\mr_+$, the set of polytopes $\mathcal{F}_\Phi$ provides us with a cellular decomposition of $X$. At this juncture, we note that one \it{can} use the corresponding cell complex to compute the cohomology $H^*(X_0(\mathfrak{n}),\mc)$ directly. More precisely, there is a retract $W$ of the space $X$, known as the \it{well-rounded retract}, from which one can obtain a cell complex using the decomposition induced by that on $X$, and then one can compute $H^*(X_0(\mathfrak{n}),\mc)$ by considering the $\Gamma_0(\mathfrak{n})$-\it{equivariant} cohomology $H^*_{\Gamma_0(\mathfrak{n})}(W,\mc)$. 

Since $\Gamma_0(\mathfrak{n})$ has finite index in $\mathrm{GL}_2(\mathcal{O})$, there are only finitely many cells of a given dimension up to $\Gamma_0(\mathfrak{n})$-equivalence, and thus one can compute the cohomology using standard techniques for finite cell complexes (as is done in \textbf{[GY13]}, for example).

We prefer to work completely with the sharbly complex, but still make use of the cellular decomposition of $X$ that we have just described. Note that any $k$-sharbly $\mathbf{u} = [u_1,\ldots,u_{k+2}]$ defines a convex polytope $\mathcal{P}(\mathbf{u})$ in $\mathcal{V}$ by taking the convex hull of the points $q(u_i)$. We define a $k$-sharbly $\mathbf{u}$ to be
\begin{itemize} 
\item\it{reduced} if $\mathcal{P}(\mathbf{u})$ is contained in $\mathcal{F}_\Phi$ for some perfect form $\Phi$; and
\item\it{totally reduced} if $\mathcal{P}(\mathbf{u})$ defines a $(k+1)$-dimensional face of $\mathcal{F}_\Phi$ for some perfect form $\Phi$.
\end{itemize}

Since there are only finitely many perfect pyramids under the action of $\Gamma_0(\mathfrak{n})$, it follows that there are only finitely many reduced and totally reduced sharblies up to $\Gamma_0(\mathfrak{n})$-equivalence. In \textbf{[GHY13]}, the authors use the subcomplex of reduced sharblies to compute the homology of the sharbly complex. We choose to work instead with the smaller subcomplex of totally reduced sharblies, which our data suggests also computes the homology of $(\mathscr{S}_*)_{\Gamma_0(\mathfrak{n})}$ in practice.

It soon becomes apparent that the subcomplex of totally reduced sharblies is not preserved by the action of the Hecke operators $T_v$ defined previously. Thus we would like to find a means for rewriting an arbitrary $k$-sharbly chain in terms of totally reduced sharblies, in order to be able to compute this action. While this is not possible in general, it turns out that if $k \in \{0,1\}$ then, given a $k$-sharbly \it{cycle} representing a class in $H_k((\mathscr{S}_*)_{\Gamma_0(\mathfrak{n})},\mc)$, we can in practice find another representative of the same class whose support consists entirely of totally reduced sharblies, using an algorithm described in \textbf{[Gun00]}, which we shall now briefly describe.

To begin with, suppose we are given a generic $0$-sharbly $\mathbf{u} = [u_1,u_2]$, and consider the point $\mathcal{B}(\mathbf{u}) = \tfrac{1}{2}(q(u_1) + q(u_2))$. Since $\overline{\mathcal{C}}$ is a convex cone, $\mathcal{B}(\mathbf{u})$ lies inside it, and in particular is contained in some perfect pyramid, corresponding to the perfect form $\Phi$, say. Define the \it{size} $N(\mathbf{u})$ of $\mathbf{u}$ to be the inner product $\langle \Phi, \mathcal{B}(\mathbf{u})\rangle$, which is well-defined if we impose the condition that the vectors $u_i \in \mathcal{O}^2$ are chosen such that no point of $\Xi$ lies on the line segment joining $q(u_i)$ and the origin. It is not too difficult to see that for each such line there exists a unique point of $\Xi$ satisfying this requirement, and by the relations defining the sharbly complex we can assume that this condition is satisfied by all sharblies.

By definition, $N(\mathbf{u}) \geq 1$, with equality if, and only if, the points $q(u_1)$ and $q(u_2)$ are minimal vectors of the form $\Phi$, which occurs if, and only if, $\mathbf{u}$ is reduced. The size of a $0$-sharbly therefore seems an appropriate measure for how far it is from being reduced.

Now, suppose that $\mathbf{u}$ is not reduced, and note that for any point $x \in \mathcal{O}^2$ we have \begin{equation*}\mathbf{u} + \partial [u_2,u_1,x] = [u_1,x] + [x,u_2],\end{equation*} and thus $\mathbf{u}$ is homologous to $[u_1,x] + [x,u_2]$. We call a point $x \in \mathcal{O}^2$ a \it{reducing point} for $\mathbf{u}$ if \begin{equation*}\mathrm{Max}\{N([u_1,x]),N([x,u_2])\} < N(\mathbf{u}).\end{equation*} One can show that, if $\Phi$ denotes the perfect form corresponding to the perfect pyramid containing the point $\mathcal{B}(\mathbf{u})$, then there exists a reducing point $x$ for $\mathbf{u}$ such that $q(x)$ is a minimal vector for $\Phi$.

One can visualize this process by means of the following diagram. Here a $0$-sharbly is represented by a straight line, with shorter lines representing sharblies of smaller size: 

\begin{center}
\begin{tikzpicture}[scale=1]
\draw (0,0)--(1.6,0)--(3.2,0);
\node [below] at (0,0) {$u_1$};
\node [below] at (1.6,0) {$x$};
\node [below] at (3.2,0) {$u_2$};
\draw [fill] (0,0) circle [radius=0.04cm];
\draw [fill] (1.6,0) circle [radius=0.04cm];
\draw [fill] (3.2,0) circle [radius=0.04cm];
\end{tikzpicture}
\end{center}

If $\mathbf{u}$ is reduced, but not totally reduced, we take a slightly different approach. One can define an alternative notion of the size of $\mathbf{u}$ (as in \textbf{[GHY13]}) which we shall denote by $n(\mathbf{u})$, by setting $n(\mathbf{u}) = \mathrm{Norm}_{F/\mq}\left(\mathrm{det}(u_1|u_2)\right)$. There appears to be a correlation between the values $N(\mathbf{u})$ and $n(\mathbf{u})$, with small values of $N(\mathbf{u})$ typically corresponding to small values of $n(\mathbf{u})$ - for example, over the field $\mq(\zeta_{12})$, we observe that a non-degenerate $0$-sharbly $\mathbf{u}$ is reduced if, and only if, $n(\mathbf{u}) \in \{1,4,9\}$, and is totally reduced if, and only if, $n(\mathbf{u}) = 1$. We therefore call a point $x \in \mathcal{O}^2$ a \it{reducing point} for a reduced $0$-sharbly $\mathbf{u}$ if $\mathrm{Max}\{n([u_1,x]),n([x,u_2])\} < n(\mathbf{u})$. Over the field $\mq(\zeta_{12})$, one can always find such a reducing point.

We now consider the case of $1$-sharblies, which is considerably more involved. Our original notion of size can be extended to a $1$-sharbly $\mathbf{u} = [u_1,u_2,u_3]$ by defining $N(\mathbf{u})$ to be $\langle \Phi, \mathcal{B}(\mathbf{u})\rangle$, where $\mathcal{B}(\mathbf{u}) = \tfrac{1}{3}(q(u_1)+q(u_2)+q(u_3))$, and $\Phi$ once again denotes the perfect form corresponding to the perfect pyramid within which $\mathcal{B}(\mathbf{u})$ lies. Once again, $\mathbf{u}$ is reduced if, and only if, $N(\mathbf{u}) = 1$. (We note that one can similarly extend the definition of $n(\mathbf{u})$ to cater to $1$-sharblies by setting $n(\mathbf{u}) = \mathrm{Max}\{n([u_i,u_j])\}$).

Now, suppose we are given a generic $1$-sharbly $\mathbf{u} = [u_1,u_2,u_3]$, none of whose edges are totally reduced, and choose reducing points $x_1$, $x_2$ and $x_3$ for the edges $[u_2,u_3]$, $[u_3,u_1]$ and $[u_1,u_2]$ respectively. Define a $2$-sharbly chain \begin{equation*}\nu := [u_1,u_2,u_3,x_1] + [u_3,u_1,x_1,x_2] + [u_1,u_2,x_1,x_3] + [u_1,x_1,x_2,x_3],\end{equation*} so that \begin{equation*}\mathbf{u} + \partial \nu = [u_1,x_3,x_2] + [u_2,x_1,x_3] + [u_3,x_2,x_1] + [x_1,x_2,x_3] + \underset{\sigma \in A_3}{\displaystyle\sum} ~[u_{\sigma(1)},u_{\sigma(2)},x_{\sigma(3)}].\end{equation*}

By definition, $\mathbf{u} + \partial \nu$ is homologous to $\mathbf{u}$. If we were able to neglect the final terms $[u_{\sigma(1)},u_{\sigma(2)},x_{\sigma(3)}]$, then we could represent the process of replacing $\mathbf{u}$ with $\mathbf{u}+\partial \nu$ by the following diagram, with $1$-sharblies represented by triangles, where a smaller triangle corresponds to a sharbly of smaller size (we note that, unlike in the case of $0$-sharblies, we cannot prove whether these sharblies are indeed smaller than $\mathbf{u}$, although in practice they seem to be):

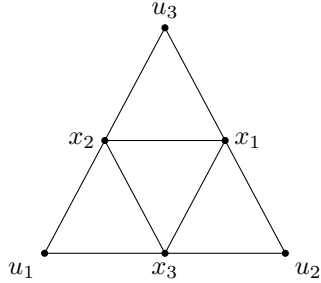
\begin{center}
\begin{tikzpicture}[scale=1]
\draw (0,0)--(1.6,3)--(3.2,0)--(0,0);
\draw (0.8,1.5)--(1.6,0)--(2.4,1.5)--(0.8,1.5);
\node [below left] at (0,0) {$u_1$};
\node [below right] at (3.2,0) {$u_2$};
\node [above] at (1.6,3) {$u_3$};
\node [left] at (0.8,1.5) {$x_2$};
\node [right] at (2.4,1.5) {$x_1$};
\node [below] at (1.6,0) {$x_3$};
\draw [fill] (0,0) circle [radius=0.04cm];
\draw [fill] (1.6,3) circle [radius=0.04cm];
\draw [fill] (3.2,0) circle [radius=0.04cm];
\draw [fill] (0.8,1.5) circle [radius=0.04cm];
\draw [fill] (1.6,0) circle [radius=0.04cm];
\draw [fill] (2.4,1.5) circle [radius=0.04cm]; 
\end{tikzpicture}
\end{center}

In general, we cannot neglect the terms $[u_{\sigma(1)},u_{\sigma(2)},x_{\sigma(3)}]$. However, if $\mathbf{u}$ lies in the support of some $1$-sharbly cycle $\xi$, then through careful choice of reducing points we can ensure that, when we perform this procedure over the entire chain, all such terms vanish. To provide a brief justification of this statement, let $\mathbf{u} \in \mathrm{Supp}(\xi)$, and consider the single edge $[u_1,u_2]$ of $\mathbf{u}$. Since the boundary of $\xi$ vanishes in $(\mathscr{S}_*)_{\Gamma_0(\mathfrak{n})}$, one can show that there must exist some $\mathbf{v} = [v_1,v_2,v_3] \in \mathrm{Supp}(\xi)$ and some edge $[v_i,v_j]$ of $\mathbf{v}$ such that \begin{equation*}[u_1,u_2] + [v_i,v_j] = 0~(\mathrm{mod}~\Gamma_0(\mathfrak{n})),\end{equation*} or equivalently that \begin{equation*}[v_i,v_j] = \gamma [u_2,u_1]\end{equation*} for some $\gamma \in \Gamma_0(\mathfrak{n})$. Let $x$ be a reducing point for $[u_1,u_2]$, and assign to $[v_i,v_j]$ the reducing point $\gamma  x$. Then \begin{equation*}[u_1,u_2,x] + [v_1,v_2,\gamma x] = [u_1,u_2,x] + \gamma [u_2,u_1,x] = 0~(\mathrm{mod}~\Gamma_0(\mathfrak{n})).\end{equation*}

\newpage Thus, as long as we can choose our reducing points in this manner, the term \begin{equation*}\underset{\mathbf{u} \in \mathrm{Supp}(\xi)}{\displaystyle\sum}~\underset{\sigma \in A_3}{\displaystyle\sum}~[u_{\sigma(1)},u_{\sigma(2)},x_{\sigma(3)}]\end{equation*} vanishes in $(\mathscr{S}_1)_{\Gamma_0(\mathfrak{n})}$ for any $1$-sharbly cycle $\xi$, so we need a means of ensuring that these points are chosen appropriately. We do this as follows: for any $0$-sharbly $\mathbf{v}$, define a \it{lift} of $\mathbf{v}$ to be a matrix $M_\mathbf{v} \in M_2(\mathcal{O})$ such that, if $m_1$ and $m_2$ denote the columns of $M_\mathbf{v}$, then $[m_1,m_2] = \pm \mathbf{v}$, subject to the relations defining the sharbly complex. Now, to each $1$-sharbly in $\mathrm{Supp}(\xi)$, we assign three lift matrices, one for each edge, such that if $\mathbf{v}$ and $\mathbf{w}$ are two edges, satisfying $\mathbf{v} = -\gamma \mathbf{w}$ for some $\gamma \in \Gamma_0(\mathfrak{n})$, then $M_\mathbf{v} = \gamma M_\mathbf{w}$.

Now, suppose we wish to choose a reducing point for an edge $\mathbf{v}$ of $\mathbf{u} \in \mathrm{Supp}(\xi)$. If $M_\mathbf{v} = \gamma M_\mathbf{w}$, where $\gamma \in \Gamma_0(\mathfrak{n})$ and $\mathbf{w}$ is an edge for which we have already chosen a reducing point $x_\mathbf{w}$, then we assign to $\mathbf{v}$ the reducing point $\gamma x_\mathbf{w}$. Otherwise, we choose an arbitrary reducing point $x_\mathbf{v}$ for $\mathbf{v}$, and record both $M_\mathbf{v}$ and $x_\mathbf{v}$. In particular, if $\mathbf{v} = -\gamma \mathbf{w}$, then $M_\mathbf{v} = \gamma M_\mathbf{w}$, and so $x_\mathbf{v} = \gamma x_\mathbf{w}$, as we have required.

In practice, repeated applications of this procedure (with some minor adjustments to cater for certain exceptional cases - for which we refer the reader to the expositions in \textbf{[Gun11]} or \textbf{[Jon14]}) eventually produce a sharbly chain supported entirely by totally reduced $1$-sharblies. In particular, given a basis for the homology group $H_1((\mathscr{S}_*)_{\Gamma_0(\mathfrak{n})},\mc)$ and a non-archimedean place $v$ of $F$, we can compute the action of the Hecke operator $T_v$ on the sharbly homology.

\newpage\specialsection{\textbf{Results}}

The following pages give details of cuspidal Hecke eigenclasses defined over the field $F = \mq(\zeta_{12})$, for which we investigated the cohomology $H^5(X_0(\mathfrak{n}),\mc)$ for all levels $\mathfrak{n}$ with norm at most $5500$. In order to detect non-trivial cuspidal cohomology, we required data regarding the rank of the corresponding Eisenstein cohomology. This is provided in the following table, which collates heuristic data from \textbf{[GHY13]} (to determine these ranks, one observes that the Hecke operator $T_v$, where $v$ does not divide the level $\mathfrak{n}$, acts on the Eisenstein subspace via multiplication by $\mathrm{Norm}_{F/\mq}(v)+1$):

\begin{table}[h]\centering\footnotesize\begin{tabular}{|c||c|c|c|c|c|c|c|}\hline Factorisation Type&$\mathfrak{p}$&$\mathfrak{p}^2$&$\mathfrak{p}^3$&$\mathfrak{p}^4$&$\mathfrak{p}^5$&$\mathfrak{p}^6$&$\mathfrak{pq}$\\\hline$\mathrm{dim}~ H^5_{\mathrm{Eis}}(X_0(\mathfrak{n}),\mc)$&$3$&$5$&$7$&$9$&$11$&$13$&$7$\\\hhline{|=#=|=|=|=|=|=|=|}
Factorisation Type&$\mathfrak{p}^2\mathfrak{q}$&$\mathfrak{p}^3\mathfrak{q}$&$\mathfrak{p}^4\mathfrak{q}$&$\mathfrak{p}^2\mathfrak{q}^2$&$\mathfrak{p}^3\mathfrak{q}^2$&$\mathfrak{pqr}$&$\mathfrak{p}^2\mathfrak{qr}$\\\hline
$\mathrm{dim}~ H^5_{\mathrm{Eis}}(X_0(\mathfrak{n}),\mc)$&$11$&$15$&$19$&$17$&$23$&$15$&$23$\\\hline
\end{tabular}\end{table}

Up to Galois conjugation, there are $544$ levels whose norm lies within the studied range. We discovered non-Eisenstein cohomology at $55$ of these levels, with a total of $99$ non-Eisenstein Hecke eigenclasses spread across these levels. Table 5.1 below lists a set of generators for the levels studied, together with their factorization type and the discrepancy $d$ between the rank of $H_1((\mathscr{S}_*)_{\Gamma_0(\mathfrak{n})},\mc)$ and the expected rank of the Eisenstein cohomology (here $t = \zeta_{12}$ denotes a primitive twelfth root of unity):

\begin{table}[h]\centering\scriptsize\begin{tabular}{|c|c|c|c||c|c|c|c|}\hline Level&Generator&Type&$d$&Level&Generator&Type&$d$\\\hhline{|=|=|=|=#=|=|=|=|}
$169$&$2t^3-3t^2-3t+2$&$\mathfrak{pq}$&$1$&$3721a$&$7t^3-6t^2-t-1$&$\mathfrak{pq}$&$3$\\\hline
$441$&$5t^2-1$&$\mathfrak{pq}$&$1$&$3721b$&$6t^2-5t-6$&$\mathfrak{pq}$&$1$\\\hline
$484$&$t^3+4t^2-4t-1$&$\mathfrak{pq}$&$1$&$3844$&$5t^3-t^2+t+6$&$\mathfrak{pq}$&$1$\\\hline
$576$&$2t^3+2t^2+2t-4$&$\mathfrak{p}^3\mathfrak{q}$&$1$&$3969$&$9t^2-6$&$\mathfrak{p}^2\mathfrak{q}$&$2$\\\hline
$625$&$5$&$\mathfrak{pq}$&$2$&$4033a$&$-8t^3+9t-9$&$\mathfrak{pq}$&$1$\\\hline
$676$&$3t^3-t^2+3t$&$\mathfrak{pqr}$&$2$&$4033b$&$-11t^3+6t^2+5t-9$&$\mathfrak{pq}$&$1$\\\hline
$1089$&$-t^3+2t-6$&$\mathfrak{pq}$&$2$&$4057$&$6t^3+2t^2-9t-2$&$\mathfrak{p}$&$1$\\\hline
$1156$&$3t^2+5t-3$&$\mathfrak{pq}$&$1$&$4069$&$-7t^3-6t^2+6t+2$&$\mathfrak{pq}$&$1$\\\hline
$1369$&$2t^3+2t^2+3t-5$&$\mathfrak{pq}$&$2$&$4096$&$8$&$\mathfrak{p}^6$&$1$\\\hline
$1521$&$4t^3+4t^2-5t+1$&$\mathfrak{pqr}$&$2$&$4225a$&$-5t^3+3t^2+9t-3$&$\mathfrak{pqr}$&$2$\\\hline
$1764$&$t^3+4t^2+4t-5$&$\mathfrak{pqr}$&$2$&$4225b$&$-9t^3+3t^2+6t-1$&$\mathfrak{p}^2\mathfrak{q}$&$1$\\\hline
$1936$&$4t^3-4t^2-6t-2$&$\mathfrak{p}^2\mathfrak{q}$&$2$&$4225c$&$-4t^3+7$&$\mathfrak{pqr}$&$1$\\\hline
$2041$&$-t^3+6t^2-t-7$&$\mathfrak{pq}$&$1$&$4356$&$5t^3+3t^2+5t$&$\mathfrak{pqr}$&$6$\\\hline
$2116$&$5t^3-5t^2+t+6$&$\mathfrak{pq}$&$2$&$4516$&$-4t^3-3t^2+9t+1$&$\mathfrak{pq}$&$1$\\\hline
$2197a$&$t^3-2t^2+3t+7$&$\mathfrak{pqr}$&$2$&$4624$&$-8t^3-2$&$\mathfrak{p}^2\mathfrak{q}$&$2$\\\hline
$2197b$&$t^3+2t^2-7t-2$&$\mathfrak{p}^2\mathfrak{q}$&$2$&$4672$&$8t^3+6t^2-6t-2$&$\mathfrak{p}^3\mathfrak{q}$&$1$\\\hline
$2209$&$4t^3-8t-1$&$\mathfrak{p}$&$1$&$4761$&$-7t^3+5t^2+2t+2$&$\mathfrak{pq}$&$3$\\\hline
$2257$&$-2t^3+6t^2+5t+1$&$\mathfrak{pq}$&$1$&$4852$&$-4t^3+7t^2+3t+1$&$\mathfrak{pq}$&$1$\\\hline
$2304$&$8t^3-4t$&$\mathfrak{p}^4\mathfrak{q}$&$2$&$5041$&$-8t^3+3t^2+3t-8$&$\mathfrak{p}$&$2$\\\hline
$2401$&$7$&$\mathfrak{pq}$&$3$&$5184$&$-6t^2+6t+6$&$\mathfrak{p}^3\mathfrak{q}^2$&$2$\\\hline
$2452$&$-7t^3+t^2+t+2$&$\mathfrak{pq}$&$1$&$5317$&$-7t^3+3t^2-2t-4$&$\mathfrak{pq}$&$1$\\\hline
$2500a$&$-t^3-7t^2+t$&$\mathfrak{p}^2\mathfrak{q}$&$1$&$5329a$&$t^3+5t^2+3t-9$&$\mathfrak{p}^2$&$2$\\\hline
$2500b$&$5t^2+5t-5$&$\mathfrak{pqr}$&$4$&$5329b$&$3t^3-8t^2-3t$&$\mathfrak{pq}$&$2$\\\hline
$2704$&$-2t^3-6t^2+6t+2$&$\mathfrak{p}^2\mathfrak{qr}$&$4$&$5329c$&$3t^3-6t-10$&$\mathfrak{pq}$&$4$\\\hline
$2916$&$3t^3-3t^2+3t+6$&$\mathfrak{p}^3\mathfrak{q}$&$2$&$5329d$&$8t^3-9t$&$\mathfrak{pq}$&$1$\\\hline
$2977$&$4t^3+2t^2-9t+2$&$\mathfrak{pq}$&$1$&$5473$&$-9t-8$&$\mathfrak{pq}$&$1$\\\hline
$3328$&$4t^3+8t^2-4t-4$&$\mathfrak{p}^4\mathfrak{q}$&$1$&$5476$&$-5t^2-t-5$&$\mathfrak{pqr}$&$5$\\\hline
$3481$&$5t^3-5t^2-6t-1$&$\mathfrak{p}$&$2$&&&&\\\hline
\end{tabular}\caption*{Table 5.1: Levels with non-Eisenstein cohomology classes}\end{table}

\newpage Of the $99$ non-Eisenstein Hecke eigenclasses we detected:
\begin{itemize}
\item $68$ admitted rational eigenvalues;
\item $18$ admitted eigenvalues lying in a quadratic extension of $\mq$;
\item $9$ admitted eigenvalues lying in a cubic extension of $\mq$; and
\item $4$ admitted eigenvalues lying in a quartic extension of $\mq$.
\end{itemize}

Of the $68$ eigenclasses which admitted rational Hecke eigenvalues:
\begin{itemize}
\item $31$ had eigenvalues matching an eigenclass appearing at a lower level;
\item $15$ had eigenvalues matching those expected from the base change of an automorphic form defined over a quadratic subfield of $F$;
\item $2$ had eigenvalues matching those from the Eisenstein cohomology, up to sign; and
\item $20$ classes could not be attributed to any of these phenomena, and we were able to find elliptic curves defined over $F$ whose local data matched the eigenvalue data for each of these classes. We list these classes in Table 5.2 below, together with their Hecke eigenvalues for a number of primes of small norm, while the corresponding elliptic curves are listed in Table 5.10.
\end{itemize}

\begin{table}[h]\centering\scriptsize\begin{tabular}{|c||c||c||c|c|c|c||c|c|}\hline Class&$\mathfrak{p}_2$&$\mathfrak{p}_3$&$\mathfrak{p}_{13,1}$&$\mathfrak{p}_{13,2}$&$\mathfrak{p}_{13,3}$&$\mathfrak{p}_{13,4}$&$\mathfrak{p}_{5,1}$&$\mathfrak{p}_{5,2}$\\\hhline{|=#=#=#=|=|=|=#=|=|}
$441$&$0$&$*$&$-6$&$4$&$4$&$-6$&$-4$&$-4$\\\hline
$1156$&$*$&$0$&$4$&$4$&$-6$&$-6$&$6$&$6$\\\hline
$2041$&$2$&$-2$&$2$&$2$&$*$&$-4$&$-4$&$-10$\\\hline
$2257$&$-3$&$-4$&$-1$&$1$&$-6$&$-3$&$1$&$-8$\\\hline
$2452$&$*$&$1$&$-4$&$-4$&$-4$&$5$&$-1$&$8$\\\hline
$2500a$&$*$&$0$&$4$&$4$&$-1$&$-1$&$1$&$*$\\\hline
$2977$&$2$&$4$&$2$&$-4$&$*$&$-4$&$2$&$-10$\\\hline
$3328$&$*$&$2$&$-2$&$-2$&$6$&$*$&$2$&$-6$\\\hline
$3721b$&$2$&$-2$&$-4$&$-4$&$2$&$2$&$8$&$8$\\\hline
$3844$&$*$&$-5$&$-1$&$-6$&$-6$&$-1$&$1$&$1$\\\hline
$4033a$&$-1$&$4$&$2$&$-4$&$2$&$2$&$2$&$-4$\\\hline
$4033b$&$2$&$-2$&$-4$&$2$&$2$&$2$&$-4$&$2$\\\hline
$4057$&$-3$&$-2$&$-4$&$-1$&$-4$&$1$&$-5$&$-2$\\\hline
$4069$&$-3$&$-4$&$-3$&$1$&$*$&$-5$&$7$&$1$\\\hline
$4225b$&$-2$&$-2$&$-4$&$*$&$-2$&$-6$&$*$&$4$\\\hline
$4516$&$*$&$5$&$4$&$-1$&$-6$&$4$&$-4$&$6$\\\hline
$4672$&$*$&$2$&$-2$&$-2$&$-2$&$6$&$-6$&$2$\\\hline
$4852$&$*$&$-3$&$-1$&$-7$&$-2$&$-4$&$3$&$-8$\\\hline
$5317$&$-3$&$2$&$-2$&$6$&$-2$&$*$&$2$&$2$\\\hline
$5473$&$-1$&$-2$&$2$&$2$&$-4$&$*$&$2$&$8$\\\hline
\end{tabular}\caption*{Table 5.2: Rational Hecke eigenclasses over $F$}\end{table}

Table 5.3 below lists generators for the prime ideals of norm up to $25$:

\begin{table}[h]\centering\small\begin{tabular}{|c|c||c|c|}\hline$\mathfrak{p}$&Generator&$\mathfrak{p}$&Generator\\\hhline{|=|=#=|=|}
$\mathfrak{p}_2$&$-t^2+t+1$&$\mathfrak{p}_{13,3}$&$-t^3-t+1$\\\cline{1-2}
$\mathfrak{p}_3$&$t^2+1$&$\mathfrak{p}_{13,4}$&$t^3+t^2+1$\\\hline
$\mathfrak{p}_{13,1}$&$-t^3+t^2+1$&$\mathfrak{p}_{5,1}$&$2t^2-t-2$\\
$\mathfrak{p}_{13,2}$&$t^3+t+1$&$\mathfrak{p}_{5,2}$&$t^3-2t^2-t$\\\hline
\end{tabular}\caption*{Table 5.3: Generators for prime ideals of $F$ of small norm}\end{table}

The classes listed in Table 5.4 below are ``old'', in the sense that the Hecke eigenvalues of this class match those of a non-Eisenstein class appearing at a level $\mathfrak{d}$ dividing $\mathfrak{n}$. Lower case Roman numerals are used to denote each eigenclass.

\begin{table}[h]\tiny\begin{tabular}{|c||c||c||c|c|c|c||c|c||c|}\hline Class&$\mathfrak{p}_2$&$\mathfrak{p}_3$&$\mathfrak{p}_{13,1}$&$\mathfrak{p}_{13,2}$&$\mathfrak{p}_{13,3}$&$\mathfrak{p}_{13,4}$&$\mathfrak{p}_{5,1}$&$\mathfrak{p}_{5,2}$&Original Class\\\hhline{|=#=#=#=|=|=|=#=|=#=|}
$676$ (\textbf{i-ii})&$*$&$-4$&$0$&$*$&$0$&$*$&$-2$&$-2$&$169$\\\hline
$1521$ (\textbf{i-ii})&$-2$&$*$&$0$&$*$&$0$&$*$&$-2$&$-2$&$169$\\\hline
$1764$ (\textbf{i-ii})&$*$&$*$&$-6$&$4$&$4$&$-6$&$-4$&$-4$&$441$\\\hline
$1936$ (\textbf{i-ii})&$*$&$-5$&$-1$&$-1$&$-1$&$-1$&$-4$&$-4$&$484$\\\hline
$2197a$ (\textbf{i-ii})&$-2$&$-4$&$0$&$*$&$*$&$*$&$-2$&$-2$&$169$\\\hline
$2197b$ (\textbf{i-ii})&$-2$&$-4$&$0$&$*$&$0$&$*$&$-2$&$-2$&$169$\\\hline
$2304$ (\textbf{i-ii})&$*$&$*$&$-2$&$-2$&$-2$&$-2$&$-6$&$-6$&$576$\\\hline
$2704$ (\textbf{ii-iv})&$*$&$-4$&$0$&$*$&$0$&$*$&$-2$&$-2$&$169$\\\hline
$3969$ (\textbf{i-ii})&$0$&$*$&$-6$&$4$&$4$&$-6$&$-4$&$-4$&$441$\\\hline
$4225a$ (\textbf{i-ii})&$-2$&$-4$&$0$&$*$&$0$&$*$&$-2$&$*$&$169$\\\hline
$4356$ (\textbf{i-ii})&$*$&$*$&$-2$&$6$&$-2$&$6$&$-6$&$-6$&$1089$ (I)\\
$4356$ (\textbf{iii-iv})&$*$&$*$&$4$&$-6$&$4$&$-6$&$6$&$6$&$1089$ (II)\\
$4356$ (\textbf{v-vi})&$*$&$*$&$-1$&$-1$&$-1$&$-1$&$-4$&$-4$&$484$\\\hline
$4624$ (\textbf{i-ii})&$*$&$0$&$4$&$4$&$-6$&$-6$&$6$&$6$&$1156$\\\hline
$5184$ (\textbf{i-ii})&$*$&$*$&$-2$&$-2$&$-2$&$-2$&$-6$&$-6$&$576$\\\hline
\end{tabular}\caption*{Table 5.4: ``Old'' cohomology classes}\end{table}

The classes listed in Table 5.5 below correspond to the base change of an automorphic representation $\pi'$ defined over a subfield $F'$ of $F$, such that the Hecke eigenvalues $a_\mathfrak{q}(\pi')$ are rational. If $\pi$ is a base change of $\pi'$, and $\mathfrak{p}$ is a prime of $F$ lying above a prime $\mathfrak{q}$ of $F'$, then \begin{equation*}a_\mathfrak{p}(\pi) = \left\{\begin{array}{cc}a_\mathfrak{q}(\pi');&~\mathrm{if~}\mathfrak{q}~\mathrm{splits~in~}F,\\a_\mathfrak{q}(\pi')^2-2\mathrm{Norm}_{F'/\mq}(\mathfrak{q});&~\mathrm{if~}\mathfrak{q}~\mathrm{is~inert~in~}F.\end{array}\right.\end{equation*}
For each of these classes, we were able to find an elliptic curve defined over the corresponding subfield whose local data matched these eigenvalues (listed in Table 5.11). 

\begin{table}[h]\centering\tiny\begin{tabular}{|c||c||c||c|c|c|c||c|c||c|}\hline Class&$\mathfrak{p}_2$&$\mathfrak{p}_3$&$\mathfrak{p}_{13,1}$&$\mathfrak{p}_{13,2}$&$\mathfrak{p}_{13,3}$&$\mathfrak{p}_{13,4}$&$\mathfrak{p}_{5,1}$&$\mathfrak{p}_{5,2}$&Base Field\\\hhline{|=#=#=#=|=|=|=#=|=#=|}
$484$&$*$&$-5$&$-1$&$-1$&$-1$&$-1$&$-4$&$-4$&$\mq(\sqrt{3})$\\\hline
$576$&$*$&$*$&$-2$&$-2$&$-2$&$-2$&$-6$&$-6$&$\mq(\sqrt{3})$\\\hline
$1089$ (\textbf{i})&$-3$&$*$&$-2$&$6$&$-2$&$6$&$-6$&$-6$&$\mq(\sqrt{3})$\\
$1089$ (\textbf{ii})&$0$&$*$&$4$&$-6$&$4$&$-6$&$6$&$6$&$\mq(\sqrt{3})$\\\hline
$2209$&$-3$&$-2$&$-6$&$0$&$-6$&$0$&$-6$&$-6$&$\mq(\sqrt{3})$\\\hline
$2704$ (\textbf{i})&$*$&$-2$&$2$&$*$&$2$&$*$&$2$&$2$&$\mq(\sqrt{3})$\\\hline
$2916$ (\textbf{i})&$*$&$*$&$5$&$-4$&$5$&$-4$&$-1$&$-1$&$\mq(\sqrt{3})$\\
$2916$ (\textbf{ii})&$*$&$*$&$-4$&$5$&$-4$&$5$&$-1$&$-1$&$\mq(\sqrt{3})$\\\hline
$4225c$&$-4$&$-2$&$*$&$*$&$-4$&$-4$&$*$&$-10$&$\mq(\sqrt{-1})$\\\hline
$5041$ (\textbf{i})&$0$&$5$&$-6$&$-1$&$-6$&$-1$&$1$&$1$&$\mq(\sqrt{3})$\\
$5041$ (\textbf{ii})&$-4$&$5$&$2$&$-1$&$2$&$-1$&$-7$&$-7$&$\mq(\sqrt{3})$\\\hline
$5329d$&$-1$&$-2$&$2$&$2$&$2$&$2$&$2$&$2$&$\mq(\sqrt{-3})$\\\hline
$5476$ (\textbf{i})&$*$&$-5$&$-4$&$-7$&$-4$&$-7$&$2$&$2$&$\mq(\sqrt{3})$\\\hline
\end{tabular}\caption*{Table 5.5: Base change from rational Hecke eigenclasses}\end{table}

\newpage In Table 5.6 we list the remaining eigenclasses which correspond to the base change of an automorphic representation $\pi'$ defined over a subfield of $F$. In each case, the Hecke eigenvalues $a_\mathfrak{p}(\pi')$ lie in a quadratic extension of $\mq$, and so there is no elliptic curve defined over the corresponding subfield of $F$ whose local data matches these eigenvalues. However, for each class, we were able to find an elliptic curve defined over $F$ whose local data matched the eigenvalues $a_\mathfrak{p}(\pi)$ (listed in Table 5.12). 

\begin{table}[h]\centering\tiny\begin{tabular}{|c||c||c||c|c|c|c||c|c||c|}\hline Class&$\mathfrak{p}_2$&$\mathfrak{p}_3$&$\mathfrak{p}_{13,1}$&$\mathfrak{p}_{13,2}$&$\mathfrak{p}_{13,3}$&$\mathfrak{p}_{13,4}$&$\mathfrak{p}_{5,1}$&$\mathfrak{p}_{5,2}$&Base Field\\\hhline{|=#=#=#=|=|=|=#=|=#=|}
$169$&$-2$&$-4$&$0$&$*$&$0$&$*$&$-2$&$-2$&$\mq(\sqrt{3})$\\\hline
$4096$&$*$&$2$&$-2$&$-2$&$-2$&$-2$&$2$&$2$&$\mq(\sqrt{3})$\\\hline
\end{tabular}\caption*{Table 5.6: Base change from non-rational Hecke eigenclasses}\end{table}

In Table 5.7 we list the remaining two eigenclasses with rational Hecke eigenvalues, which match those of the Eisenstein cohomology, up to sign. We observe that the ray class group $Cl(\mathcal{O}_F,\mathfrak{n})$ of the corresponding level admits a single non-trivial quadratic character $\chi$, and that the Hecke eigenvalues are given by \begin{equation*}a_\mathfrak{p}(\pi) = \chi(\mathfrak{p})(\mathrm{Norm}_{F/\mq}(\mathfrak{p})+1).\end{equation*} These appear to correspond to classes denoted by $H^\cdot_{\mathrm{res}}$ in \textbf{[Har87], Section 3.2.5}.

\begin{table}[h]\centering\tiny\begin{tabular}{|c||c||c||c|c|c|c||c|c||c|c|c|c|}\hline Class&$\mathfrak{p}_2$&$\mathfrak{p}_3$&$\mathfrak{p}_{13,1}$&$\mathfrak{p}_{13,2}$&$\mathfrak{p}_{13,3}$&$\mathfrak{p}_{13,4}$&$\mathfrak{p}_{5,1}$&$\mathfrak{p}_{5,2}$&$\mathfrak{p}_{37,1}$&$\mathfrak{p}_{37,2}$&$\mathfrak{p}_{37,3}$&$\mathfrak{p}_{37,4}$\\\hhline{|=#=#=#=|=|=|=#=|=#=|=|=|=|}
$5329a$ (\textbf{i})&$-5$&$-10$&$14$&$-14$&$-14$&$-14$&$26$&$-26$&$-38$&$-38$&$38$&$38$\\
$5329a$ (\textbf{ii})&$-5$&$-10$&$14$&$-14$&$-14$&$-14$&$26$&$-26$&$-38$&$-38$&$38$&$38$\\\hline
\end{tabular}\caption*{Table 5.7: Remaining rational eigenclasses}\end{table}

In Tables 5.8 and 5.9 we list the remaining eigenclasses, whose eigenvalues lie in a proper extension of $\mq$. For the classes appearing in Table 5.8, the field $\mq(a_\mathfrak{p}(\pi))$ generated by these eigenvalues is a quadratic extension of $\mq$, and we list the pair of Galois conjugate eigenvalues for each prime. For the classes appearing in Table 5.9, the field $\mq(a_\mathfrak{p}(\pi))$ is either a cubic or a quartic extension of $\mq$, and for each prime we list the polynomial whose roots are the corresponding eigenvalues.

\begin{table}[h]\centering\tiny\begin{tabular}{|c||c||c||c|c|c|c||c|}\hline Class&$\mathfrak{p}_2$&$\mathfrak{p}_3$&$\mathfrak{p}_{13,1}$&$\mathfrak{p}_{13,2}$&$\mathfrak{p}_{13,3}$&$\mathfrak{p}_{13,4}$&$\mq(a_\mathfrak{p}(\pi))$\\\hhline{|=#=#=#=|=|=|=#=|}
$625$ (\textbf{i-ii})&$\tfrac{1\pm \sqrt{17}}{2}$&$-1 \pm \sqrt{17}$&$-1 \pm \sqrt{17}$&$-1 \pm \sqrt{17}$&$-1 \pm \sqrt{17}$&$-1 \pm \sqrt{17}$&$\mq(\sqrt{17})$\\\hline
$1369$ (\textbf{i-ii})&$\tfrac{-3\pm \sqrt{17}}{2}$&$\tfrac{-5 \pm \sqrt{17}}{2}$&$\tfrac{3\pm \sqrt{17}}{2}$&$1 \pm \sqrt{17}$&$\tfrac{3\pm \sqrt{17}}{2}$&$1 \pm \sqrt{17}$&$\mq(\sqrt{17})$\\\hline
$2116$ (\textbf{i-ii})&$*$&$-2 \pm 2\sqrt{3}$&$-1 \pm 3\sqrt{3}$&$2 \pm 2\sqrt{3}$&$-1 \pm 3\sqrt{3}$&$2 \pm 2\sqrt{3}$&$\mq(\sqrt{3})$\\\hline
$2500b$ (\textbf{i-iv})&$*$&$-1 \pm \sqrt{17}$&$-1 \pm \sqrt{17}$&$-1 \pm \sqrt{17}$&$-1 \pm \sqrt{17}$&$-1 \pm \sqrt{17}$&$\mq(\sqrt{17})$\\\hline
$3481$ (\textbf{i-ii})&$\tfrac{-5 \pm \sqrt{5}}{2}$&$\tfrac{-5 \pm 3\sqrt{5}}{2}$&$-1 \pm 2\sqrt{5}$&$\tfrac{-7 \pm 3\sqrt{5}}{2}$&$-1 \pm 2\sqrt{5}$&$\tfrac{-7 \pm 3\sqrt{5}}{2}$&$\mq(\sqrt{5})$\\\hline
$5329b$ (\textbf{i-ii})&$\pm \sqrt{7}$&$\pm 2 \sqrt{7}$&$-4$&$-4$&$1 \pm \sqrt{7}$&$1 \pm \sqrt{7}$&$\mq(\sqrt{7})$\\\hline
$5476$ (\textbf{ii-v})&$*$&$\tfrac{-5 \pm \sqrt{17}}{2}$&$1 \pm \sqrt{17}$&$\tfrac{3\pm \sqrt{17}}{2}$&$1 \pm \sqrt{17}$&$\tfrac{3\pm \sqrt{17}}{2}$&$\mq(\sqrt{17})$\\\hline
\end{tabular}\caption*{Table 5.8: Eigenclasses with eigenvalues lying in a quadratic extension of $\mq$}\end{table}

\begin{table}[h]\centering\tiny\begin{tabular}{|c||c||c||c|}\hline Class&$\mathfrak{p}_2$&$\mathfrak{p}_3$&$\mq(a_\mathfrak{p}(\pi))$\\\hhline{|=#=#=#=|}
$2401$ (\textbf{i-iii})&$x^3+2x^2-11x-20$&$x^3+2x^2-32x-80$&$\mq(x^3+x^2-8x-10)$\\\hline
$3721a$ (\textbf{i-iii})&$x^3+2x^2-9x-6$&$x^3+5x^2-x-2$&$\mq(x^3-x^2-9x+12)$\\\hline
$4761$ (\textbf{i-iii})&$x^3+3x^2-4x-4$&$*$&$\mq(x^3-x^2-4x+2)$\\\hline
$5329c$ (\textbf{i-iv})&$x^4+4x^3-3x^2-16x-8$&$x^4+8x^3+6x^2-48x-64$&$\mq(x^4-9x^2-2x+2)$\\\hline
\end{tabular}\caption*{Table 5.9: Eigenclasses with eigenvalues lying in a cubic or quartic extension of $\mq$}\end{table}

In Table 5.10 (on the following page) we list the coefficients $a_i$ of the Weierstrass polynomial \begin{equation*}y^2 + a_1xy + a_3y = x^3+a_2x^2+a_4x+6\end{equation*} defining the global minimal model of an elliptic curve $E$ over $F$ whose local data matches the Hecke eigenvalues of the classes appearing in Table 5.2. These curves were found using a combination of the MAGMA routine \textbf{EllipticCurveSearch} and the ideas found in \textbf{[DGKY14]}.

In Table 5.11 (below) we list the coefficients $a_i$ of the Weierstrass polynomial defining the global minimal model of an elliptic curve $E$ over a subfield of $F$ such that the Hecke eigenvalues of the classes appearing in Table 5.5 match the local data of the base change of $E$ to $F$.

In Table 5.12 (on the next page) we list the coefficients of the Weierstrass polynomial defining the global minimal model of an elliptic curve $E$ over $F$ whose local data matches the Hecke eigenvalues of the classes appearing in Table 5.6.

\begin{table}[h]\centering\scriptsize\begin{tabular}{|c||c|c|c|c|c|}\hline Class&$a_1$&$a_2$&$a_3$&$a_4$&$a_6$\\\hhline{|=#=|=|=|=|=|}
$484$&$\sqrt{3}$&$\sqrt{3}+1$&$\sqrt{3}$&$2\sqrt{3}+2$&$\sqrt{3}+1$\\\hline
$576$&$\sqrt{3}+1$&$-\sqrt{3}+1$&$0$&$-5\sqrt{3}-6$&$3\sqrt{3}+6$\\\hline
$1089$(\textbf{i})&$1$&$-\sqrt{3}$&$0$&$1$&$0$\\\hline
$1089$(\textbf{ii})&$\sqrt{3}+1$&$-\sqrt{3}$&$1$&$5\sqrt{3}-9$&$-6\sqrt{3}+10$\\\hline
$2209$&$1$&$-\sqrt{3}$&$1$&$-\sqrt{3}-1$&$0$\\\hline
$2704$(\textbf{i})&$0$&$\sqrt{3}-1$&$0$&$2$&$2\sqrt{3}+3$\\\hline
$2916$(\textbf{i})&$1$&$-1$&$\sqrt{3}+1$&$-23\sqrt{3}-41$&$217\sqrt{3}+377$\\\hline
$2916$(\textbf{ii})&$1$&$-1$&$\sqrt{3}+1$&$22\sqrt{3}-41$&$-218\sqrt{3}+377$\\\hline
$4225c$&$\sqrt{-1}+1$&$-\sqrt{-1}$&$\sqrt{-1}$&$1$&$0$\\\hline
$5041$(\textbf{i})&$0$&$-1$&$\sqrt{3}$&$-2\sqrt{3}-4$&$3\sqrt{3}+5$\\\hline
$5041$(\textbf{ii})&$0$&$1$&$\sqrt{3}$&$\sqrt{3}+2$&$\sqrt{3}+1$\\\hline
$5329d$&$3\sqrt{-3}$&$\sqrt{-3}+7$&$\tfrac{1}{2}(\sqrt{-3}-5)$&$4\sqrt{-3}+1$&$\tfrac{1}{2}(\sqrt{-3}-3)$\\\hline
$5476$&$1$&$-\sqrt{3}+1$&$\sqrt{3}$&$-\sqrt{3}+1$&$-\sqrt{3}+1$\\\hline
\end{tabular}\caption*{Table 5.11: Elliptic curves corresponding to the classes in Table 5.5}\end{table}

\begin{sidewaystable}\centering\scriptsize\begin{tabular}{|c||c|c|c|c|c|}\hline Class&$a_1$&$a_2$&$a_3$&$a_4$&$a_6$\\\hhline{|=#=|=|=|=|=|}
$441$&$-3t^3-3t^2+3t$&$2t+2$&$5t^3-2t^2-4t-2$&$-t^3-12t^2-7t+9$&$6t^3+6t^2-9t-3$\\\hline
$1156$&$-6t^3-6t^2+3$&$-13t^3-7t^2+20t+20$&$4t^3-7t^2-17t-13$&$-165t^3-66t^2+64t+11$&$235t^3-270t^2-288t+145$\\\hline
$2041$&$-3t^3+3$&$-t^3-7t^2+2t+1$&$-t^3-9t^2-2t-1$&$-10t^3+28t^2+31t+15$&$-262t^3-262t^2+98t+219$\\\hline
$2257$&$-3t^3-3t^2$&$-8t^3-7t^2+2t+7$&$9t^3+t^2-12t-13$&$-10t^3-29t^2-28t-14$&$68t^3+35t^2-65t-71$\\\hline
$2452$&$3t^3-3t-3$&$2t^3+5t^2+3t-1$&$5t^3-12t-9$&$29t^3+17t^2-9t-17$&$41t^3+22t^2-39t-33$\\\hline
$2500a$&$-3t^3-6t^2$&$-11t^3-13t^2+13t+14$&$4t^3-3t^2-8t-3$&$-15t^3-46t^2-11t+28$&$4t^3-11t^2-29t-3$\\\hline
$2977$&$-3t^3+3$&$-t^3-3t^2+2t-1$&$-8t^2-t-6$&$-11t^3+8t^2+5$&$2t^3-29t^2-3t$\\\hline
$3328$&$-3t^3+3t^2+3t$&$2t^3-9t-1$&$5t^2-3t-7$&$-50t^3+25t^2+54t-31$&$139t^3-87t^2-115t+121$\\\hline
$3721b$&$-3t^3+3$&$t^3-7t^2+t$&$5t^3-8t^2-3t-4$&$-14t^3+23t^2-t-15$&$28t^3-21t^2-13t+21$\\\hline
$3844$&$-3t^3-3t^2$&$-8t^3-7t^2+2t+8$&$8t^3-9t-7$&$-7399t^3-3866t^2+8173t+9088$&$-182355t^3-372229t^2-280418t-56472$\\\hline
$4033a$&$-3t$&$-5t^2-t+1$&$5t^3+4t^2-7t$&$8t^3-2t^2-t-10$&$4t^3+2t^2+11t-7$\\\hline
$4033b$&$-6t^3-3t^2+3t+3$&$-t^3+6t^2+7t$&$2t^3-5t^2-12t-2$&$-24t^3-28t^2+19t+37$&$-11t^3-24t^2-26t-20$\\\hline
$4057$&$-6t^3-6t^2+3$&$-9t^3-7t^2+16t+17$&$14t^3+7t^2-16t-19$&$-53t^3-82t^2-41t-3$&$132t^3+106t^2-90t-130$\\\hline
$4069$&$3t^3+3t^2+3$&$-11t^3-5t^2+9t+6$&$8t^2+4t-2$&$-17t^3-21t^2+19t+29$&$-15t^3-5t^2+14t+20$\\\hline
$4225b$&$3t^2-3t-3$&$t^3-t^2-2t$&$-4t^3+4t^2+7t+2$&$-4t^3-2t^2+5t+16$&$-2t^3-6t^2-16t-8$\\\hline
$4516$&$3t^3+3t^2-3$&$t^3+4t^2+3t$&$t^3-3t^2-3t-1$&$-5t^3+13t^2+5t-22$&$-16t^3+5t^2+10t-15$\\\hline
$4672$&$3t^3+3$&$-9t^3+4t^2+10t+9$&$4t^3+8t^2+12t+14$&$130t^3-152t^2-30t+208$&$-514t^3-122t^2+1086t-764$\\\hline
$4852$&$-3t^2-6t-3$&$-7t^3-13t^2-10t-1$&$-10t^3-9t^2+6t+7$&$-47t^3-29t^2+42t+52$&$16t^3+31t^2+24t+4$\\\hline
$5317$&$-3t^2-3t$&$-4t^3-4t^2-3t+1$&$5t^3+6t^2-t-2$&$13t^3+13t^2-8t-15$&$-6t^3+13t+10$\\\hline
$5473$&$3t^2$&$-2t^2+3t+3$&$t^2+3t-1$&$-20t^3+13t^2+58t+43$&$190t^3+188t^2-54t-144$\\\hline
\end{tabular}\caption*{Table 5.10: Elliptic curves corresponding to the classes in Table 5.2}

\bigskip\begin{tabular}{|c||c|c|c|c|c|}\hline Class&$a_1$&$a_2$&$a_3$&$a_4$&$a_6$\\\hhline{|=#=|=|=|=|=|}
$169$&$3t^3-3t^2-3t$&$-t^3+2t^2-3t-2$&$-2t^3+2t^2-t+2$&$-6t^3+t^2+9t-3$&$5t^3-3t^2-t+2$\\\hline
$4096$&$2t+2$&$-t^2+t-1$&$2t^3+2$&$-2t^3+2t^2-2t$&$0$\\\hline
\end{tabular}\caption*{Table 5.12: Elliptic curves corresponding to the classes in Table 5.6}\end{sidewaystable}

\newpage

For each of the elliptic curves appearing in Table 5.10, we list a set of primes $\mathfrak{p}$ of $F$ which suffice to prove modularity of the curve (for a more detailed discussion of how these primes were determined, we refer the reader to our exposition in \textbf{[Jon14]}). The table below gives a list of the prime ideals of $F$ of norm at most $650$, together with a generator for each ideal (we use the convention that the prime $\mathfrak{p}_{p,i}$ lies above the rational prime $p$):

\begin{table}[h]\centering\scriptsize\begin{tabular}{|c|c||c|c||c|c|}\hline$\mathfrak{p}$&Generator&$\mathfrak{p}$&Generator&$\mathfrak{p}$&Generator\\\hhline{|=|=#=|=#=|=|}
$\mathfrak{p}_2$&$-t^2+t+1$&$\mathfrak{p}_{181,3}$&$4t^3+t^2-t-2$&$\mathfrak{p}_{397,1}$&$-4t^3-2t^2+3t-1$\\\cline{1-2}
$\mathfrak{p}_3$&$t^2+1$&$\mathfrak{p}_{181,4}$&$2t^3+t^2-t-4$&$\mathfrak{p}_{397,2}$&$t^3+3t^2+2t-4$\\\cline{1-4}
$\mathfrak{p}_{13,1}$&$-t^3+t^2+1$&$\mathfrak{p}_{193,1}$&$-t^3+t^2+4t-1$&$\mathfrak{p}_{397,3}$&$-t^3+3t^2-2t-4$\\
$\mathfrak{p}_{13,2}$&$t^3+t+1$&$\mathfrak{p}_{193,2}$&$3t^2-t-4$&$\mathfrak{p}_{397,4}$&$-t^3+3t^2+4t-2$\\\cline{5-6}
$\mathfrak{p}_{13,3}$&$-t^3-t+1$&$\mathfrak{p}_{193,3}$&$3t^2+t-4$&$\mathfrak{p}_{409,1}$&$-5t^3+t^2+t-1$\\
$\mathfrak{p}_{13,4}$&$t^3+t^2+1$&$\mathfrak{p}_{193,4}$&$t^3+t^2-4t-1$&$\mathfrak{p}_{409,2}$&$-t^3+t^2+t-5$\\\cline{1-4}
$\mathfrak{p}_{5,1}$&$2t^2-t-2$&$\mathfrak{p}_{229,1}$&$-3t^3+2t^2+3t+1$&$\mathfrak{p}_{409,3}$&$-3t^3+4t^2-t-4$\\
$\mathfrak{p}_{5,2}$&$t^3-2t^2-t$&$\mathfrak{p}_{229,2}$&$-t^3+3t^2-2t-3$&$\mathfrak{p}_{409,4}$&$4t^3-5t-1$\\\cline{1-2}\cline{5-6}
$\mathfrak{p}_{37,1}$&$2t^3+t^2-2$&$\mathfrak{p}_{229,3}$&$2t^3-3t^2-3t$&$\mathfrak{p}_{421,1}$&$3t^3+t^2-2t-5$\\
$\mathfrak{p}_{37,2}$&$t^3-2t^2-2t$&$\mathfrak{p}_{229,4}$&$-3t^3-2t^2+3t-1$&$\mathfrak{p}_{421,2}$&$2t^3+4t^2-3t-5$\\\cline{3-4}
$\mathfrak{p}_{37,3}$&$t^3+t^2-t-3$&$\mathfrak{p}_{241,1}$&$4t^3-4t-1$&$\mathfrak{p}_{421,3}$&$5t^3-3t-2$\\
$\mathfrak{p}_{37,4}$&$3t^3+t^2-t-1$&$\mathfrak{p}_{241,2}$&$-t^3+4t^2-4$&$\mathfrak{p}_{421,4}$&$-3t^3+t^2+2t-5$\\\cline{1-2}\cline{5-6}
$\mathfrak{p}_{7,1}$&$2t^3-3t$&$\mathfrak{p}_{241,3}$&$t^3-t-4$&$\mathfrak{p}_{433,1}$&$-3t^3+2t^2+t-5$\\
$\mathfrak{p}_{7,2}$&$t^3-3t$&$\mathfrak{p}_{241,4}$&$-4t^3+t^2-1$&$\mathfrak{p}_{433,2}$&$5t^3+t^2-2t-3$\\\cline{1-4}
$\mathfrak{p}_{61,1}$&$-t^3+t^2+3t-1$&$\mathfrak{p}_{277,1}$&$3t^2+2t-4$&$\mathfrak{p}_{433,3}$&$t^3-2t^2-5t$\\
$\mathfrak{p}_{61,2}$&$2t^2-t-3$&$\mathfrak{p}_{277,2}$&$2t^3-t^2-3t-4$&$\mathfrak{p}_{433,4}$&$3t^3+2t^2-t-5$\\\cline{5-6}
$\mathfrak{p}_{61,3}$&$2t^2+t-3$&$\mathfrak{p}_{277,3}$&$t^3+2t^2-4t-2$&$\mathfrak{p}_{457,1}$&$3t^2+3t-4$\\
$\mathfrak{p}_{61,4}$&$t^3+t^2-3t-1$&$\mathfrak{p}_{277,4}$&$2t^3-4t^2+t-1$&$\mathfrak{p}_{457,2}$&$-3t^3-3t^2+3t-1$\\\cline{1-4}
$\mathfrak{p}_{73,1}$&$-t^3-3t^2$&$\mathfrak{p}_{17,1}$&$4t^2-t-4$&$\mathfrak{p}_{457,3}$&$3t^3-3t^2-3t-1$\\
$\mathfrak{p}_{73,2}$&$2t^3+2t^2-3$&$\mathfrak{p}_{17,2}$&$t^3-4t^2-t$&$\mathfrak{p}_{457,4}$&$t^3+3t^2+3t-3$\\\cline{3-6}
$\mathfrak{p}_{73,3}$&$-2t^3+2t^2-3$&$\mathfrak{p}_{313,1}$&$-t^3+2t^2+3t-5$&$\mathfrak{p}_{23,1}$&$2t^3+3t^2-2t-6$\\
$\mathfrak{p}_{73,4}$&$-t^3+3t^2+t-3$&$\mathfrak{p}_{313,2}$&$t^3-4t^2-2t$&$\mathfrak{p}_{23,2}$&$-2t^3+3t^2+2t-6$\\\cline{1-2}\cline{5-6}
$\mathfrak{p}_{97,1}$&$2t^3+t^2-2t-4$&$\mathfrak{p}_{313,3}$&$3t^3+3t^2-t-5$&$\mathfrak{p}_{541,1}$&$-2t^3+2t-5$\\
$\mathfrak{p}_{97,2}$&$-2t^3+4t^2-1$&$\mathfrak{p}_{313,4}$&$t^3+2t^2-3t-5$&$\mathfrak{p}_{541,2}$&$5t^3-3t-3$\\\cline{3-4}
$\mathfrak{p}_{97,3}$&$-2t^3-4t^2+1$&$\mathfrak{p}_{337,1}$&$5t^3+t^2-2t-2$&$\mathfrak{p}_{541,3}$&$5t^3-5t-2$\\
$\mathfrak{p}_{97,4}$&$-2t^3+t^2+2t-4$&$\mathfrak{p}_{337,2}$&$-2t^3+2t^2+t-5$&$\mathfrak{p}_{541,4}$&$-5t^2-2t$\\\cline{1-2}\cline{5-6}
$\mathfrak{p}_{109,1}$&$t^3+t^2-2t-4$&$\mathfrak{p}_{337,3}$&$2t^3+2t^2-t-5$&$\mathfrak{p}_{577,1}$&$-2t^3+3t^2-2t-4$\\
$\mathfrak{p}_{109,2}$&$-2t^3-2t^2+2t-1$&$\mathfrak{p}_{337,4}$&$-t^3+3t^2+2t-5$&$\mathfrak{p}_{577,2}$&$-4t^3-2t^2+3t-2$\\\cline{3-4}
$\mathfrak{p}_{109,3}$&$2t^3-2t^2-2t-1$&$\mathfrak{p}_{349,1}$&$-2t^3+3t-5$&$\mathfrak{p}_{577,3}$&$4t^3-2t^2-3t-2$\\
$\mathfrak{p}_{109,4}$&$-t^3+t^2+2t-4$&$\mathfrak{p}_{349,2}$&$4t^3-t^2-2t-2$&$\mathfrak{p}_{577,4}$&$2t^3-4t^2-4t+1$\\\cline{1-2}\cline{5-6}
$\mathfrak{p}_{11,1}$&$-3t^3+t^2+t-3$&$\mathfrak{p}_{349,3}$&$t^3+2t^2+2t-4$&$\mathfrak{p}_{601,1}$&$5t^3-5t-1$\\
$\mathfrak{p}_{11,2}$&$t^3+2t^2-3t-3$&$\mathfrak{p}_{349,4}$&$2t^3-3t-5$&$\mathfrak{p}_{601,2}$&$-t^3+5t^2-5$\\\cline{1-4}
$\mathfrak{p}_{157,1}$&$-4t^3+2t-1$&$\mathfrak{p}_{19,1}$&$2t^3-5t$&$\mathfrak{p}_{601,3}$&$t^3-t-5$\\
$\mathfrak{p}_{157,2}$&$-t^3+2t^2-4$&$\mathfrak{p}_{19,2}$&$3t^3-5t$&$\mathfrak{p}_{601,4}$&$-5t^3+t^2-1$\\\cline{3-6}
$\mathfrak{p}_{157,3}$&$t^3+2t^2-4$&$\mathfrak{p}_{373,1}$&$4t^3+2t^2-2t-5$&$\mathfrak{p}_{613,1}$&$t^3+2t^2-5t-2$\\
$\mathfrak{p}_{157,4}$&$4t^3-2t-1$&$\mathfrak{p}_{373,2}$&$2t^3+3t^2-4t-5$&$\mathfrak{p}_{613,2}$&$4t^2-2t-5$\\\cline{1-2}
$\mathfrak{p}_{181,1}$&$-2t^3+t^2+t-4$&$\mathfrak{p}_{373,3}$&$5t^3-3t-1$&$\mathfrak{p}_{613,3}$&$4t^2+2t-5$\\
$\mathfrak{p}_{181,2}$&$t^3-2t^2-4t$&$\mathfrak{p}_{373,4}$&$3t^3+t^2-5t-1$&$\mathfrak{p}_{613,4}$&$-t^3+2t^2+5t-2$\\\hline
\end{tabular}\end{table}

\newpage\textbf{Class 441}
\begin{itemize}
\item The residual representation attached to the corresponding elliptic curve has trivial image.
\item The primes $\{\mathfrak{p}_{13,1}$, $\mathfrak{p}_{13,2}$, $\mathfrak{p}_{13,3}$, $\mathfrak{p}_{13,4}$, $\mathfrak{p}_{5,1}$, $\mathfrak{p}_{5,2}$, $\mathfrak{p}_{37,1}$, $\mathfrak{p}_{37,2}$, $\mathfrak{p}_{37,3}$, $\mathfrak{p}_{7,1}$, $\mathfrak{p}_{61,2}$, $\mathfrak{p}_{73,1}\}$ suffice to prove isomorphism of the residual representations.
\item The primes $\{\mathfrak{p}_{13,1}$, $\mathfrak{p}_{13,2}$, $\mathfrak{p}_{13,3}$, $\mathfrak{p}_{13,4}$, $\mathfrak{p}_{5,1}$, $\mathfrak{p}_{37,1}$, $\mathfrak{p}_{37,2}$, $\mathfrak{p}_{37,3}$, $\mathfrak{p}_{37,4}$, $\mathfrak{p}_{7,1}$, $\mathfrak{p}_{61,2}$, $\mathfrak{p}_{61,3}$, $\mathfrak{p}_{73,1}$, $\mathfrak{p}_{73,2}$, $\mathfrak{p}_{97,1}$, $\mathfrak{p}_{97,2}$, $\mathfrak{p}_{109,2}$, $\mathfrak{p}_{109,3}$, $\mathfrak{p}_{11,1}$, $\mathfrak{p}_{181,2}$, $\mathfrak{p}_{181,3}$, $\mathfrak{p}_{193,1}$, $\mathfrak{p}_{17,1}$, $\mathfrak{p}_{313,1}$, $\mathfrak{p}_{337,1}$, $\mathfrak{p}_{349,1}$, $\mathfrak{p}_{349,4}$, $\mathfrak{p}_{19,1}$, $\mathfrak{p}_{409,2}$, $\mathfrak{p}_{23,1}$, $\mathfrak{p}_{601,2}\}$ suffice to prove isomorphism of the full representations.
\end{itemize}

\begin{table}[h]\centering\scriptsize\begin{tabular}{|c||c|c|c|c|c|c|c|c|}\hline
$\mathfrak{p}$&$\mathfrak{p}_{13,1}$&$\mathfrak{p}_{13,2}$&$\mathfrak{p}_{13,3}$&$\mathfrak{p}_{13,4}$&$\mathfrak{p}_{5,1}$&$\mathfrak{p}_{5,2}$&$\mathfrak{p}_{37,1}$&$\mathfrak{p}_{37,2}$\\\hline
$a_\mathfrak{p}$&$-6$&$4$&$4$&$-6$&$-4$&$-4$&$-2$&$-2$\\\hhline{|=#=|=|=|=|=|=|=|=|}
$\mathfrak{p}$&$\mathfrak{p}_{37,3}$&$\mathfrak{p}_{37,4}$&$\mathfrak{p}_{7,1}$&$\mathfrak{p}_{61,2}$&$\mathfrak{p}_{61,3}$&$\mathfrak{p}_{73,1}$&$\mathfrak{p}_{73,2}$&$\mathfrak{p}_{97,1}$\\\hline
$a_\mathfrak{p}$&$-2$&$-2$&$10$&$2$&$2$&$14$&$4$&$-2$\\\hhline{|=#=|=|=|=|=|=|=|=|}
$\mathfrak{p}$&$\mathfrak{p}_{97,2}$&$\mathfrak{p}_{109,2}$&$\mathfrak{p}_{109,3}$&$\mathfrak{p}_{11,1}$&$\mathfrak{p}_{181,2}$&$\mathfrak{p}_{181,3}$&$\mathfrak{p}_{193,1}$&$\mathfrak{p}_{17,1}$\\\hline
$a_\mathfrak{p}$&$8$&$10$&$10$&$2$&$-8$&$-8$&$-26$&$-20$\\\hhline{|=#=|=|=|=|=|=|=|=|}
$\mathfrak{p}$&$\mathfrak{p}_{313,1}$&$\mathfrak{p}_{337,1}$&$\mathfrak{p}_{349,1}$&$\mathfrak{p}_{349,4}$&$\mathfrak{p}_{19,1}$&$\mathfrak{p}_{409,2}$&$\mathfrak{p}_{23,1}$&$\mathfrak{p}_{601,2}$\\\hline
$a_\mathfrak{p}$&$34$&$-22$&$-30$&$-30$&$2$&$30$&$10$&$22$\\\hline
\end{tabular}\tiny\caption*{Eigenvalues $a_\mathfrak{p}$ of the Hecke operators $T_\mathfrak{p}$ on class 441}\end{table}

\textbf{Class 1156}
\begin{itemize}
\item The residual representation attached to the corresponding elliptic curve has trivial image.
\item The prime $\mathfrak{p}_{13,1}$ suffices to prove isomorphism of the residual representations.
\item The primes $\{\mathfrak{p}_3$, $\mathfrak{p}_{13,1}$, $\mathfrak{p}_{13,2}$, $\mathfrak{p}_{13,3}$, $\mathfrak{p}_{13,4}$, $\mathfrak{p}_{5,1}$, $\mathfrak{p}_{37,1}$, $\mathfrak{p}_{37,3}$, $\mathfrak{p}_{7,1}$, $\mathfrak{p}_{73,1}$, $\mathfrak{p}_{73,3}$, $\mathfrak{p}_{97,3}$, $\mathfrak{p}_{109,2}$, $\mathfrak{p}_{109,4}$, $\mathfrak{p}_{457,1}\}$ suffice to prove isomorphism of the full representations.
\end{itemize}

\begin{table}[h]\centering\scriptsize\begin{tabular}{|c||c|c|c|c|c|c|c|c|}\hline
$\mathfrak{p}$&$\mathfrak{p}_3$&$\mathfrak{p}_{13,1}$&$\mathfrak{p}_{13,2}$&$\mathfrak{p}_{13,3}$&$\mathfrak{p}_{13,4}$&$\mathfrak{p}_{5,1}$&$\mathfrak{p}_{37,1}$&$\mathfrak{p}_{37,3}$\\\hline
$a_\mathfrak{p}$&$0$&$4$&$4$&$-6$&$-6$&$6$&$-2$&$-2$\\\hhline{|=#=|=|=|=|=|=|=|=|}
$\mathfrak{p}$&$\mathfrak{p}_{7,1}$&$\mathfrak{p}_{73,1}$&$\mathfrak{p}_{73,3}$&$\mathfrak{p}_{97,3}$&$\mathfrak{p}_{109,2}$&$\mathfrak{p}_{109,4}$&$\mathfrak{p}_{457,1}$&\\\hline
$a_\mathfrak{p}$&$-10$&$4$&$-6$&$-2$&$10$&$10$&$18$&\\\hline
\end{tabular}\caption*{Eigenvalues $a_\mathfrak{p}$ of the Hecke operators $T_\mathfrak{p}$ on class 1156}\end{table}

\textbf{Class 2041}

\begin{itemize} 
\item The residual representation attached to the corresponding elliptic curve has trivial image.
\item The prime $\mathfrak{p}_{13,1}$ suffices to prove isomorphism of the residual representations.
\item The primes $\{\mathfrak{p}_{3}$, $\mathfrak{p}_{13,1}$, $\mathfrak{p}_{13,2}$, $\mathfrak{p}_{13,4}$, $\mathfrak{p}_{5,1}$, $\mathfrak{p}_{5,2}$, $\mathfrak{p}_{37,1}$, $\mathfrak{p}_{37,2}$, $\mathfrak{p}_{37,3}$, $\mathfrak{p}_{37,4}$, $\mathfrak{p}_{61,1}$, $\mathfrak{p}_{61,2}$, $\mathfrak{p}_{61,3}$, $\mathfrak{p}_{73,1}$, $\mathfrak{p}_{73,3}$, $\mathfrak{p}_{73,4}$, $\mathfrak{p}_{97,1}$, $\mathfrak{p}_{97,2}$, $\mathfrak{p}_{97,4}$, $\mathfrak{p}_{109,1}$, $\mathfrak{p}_{109,2}$, $\mathfrak{p}_{109,4}$, $\mathfrak{p}_{181,1}$, $\mathfrak{p}_{193,1}$, $\mathfrak{p}_{229,3}$, $\mathfrak{p}_{17,1}$, $\mathfrak{p}_{313,1}$, $\mathfrak{p}_{313,4}$, $\mathfrak{p}_{373,1}$, $\mathfrak{p}_{409,1}$, $\mathfrak{p}_{1321,1}\}$, where $\mathfrak{p}_{1321,1}$ is generated by the element $-3t^3-8t^2-t+3$, suffice to prove isomorphism of the full representations.
\end{itemize}

\newpage
\begin{table}[h]\centering\scriptsize\begin{tabular}{|c||c|c|c|c|c|c|c|c|}\hline
$\mathfrak{p}$&$\mathfrak{p}_{3}$&$\mathfrak{p}_{13,1}$&$\mathfrak{p}_{13,2}$&$\mathfrak{p}_{13,4}$&$\mathfrak{p}_{5,1}$&$\mathfrak{p}_{5,2}$&$\mathfrak{p}_{37,1}$&$\mathfrak{p}_{37,2}$\\\hline
$a_\mathfrak{p}$&$-2$&$2$&$2$&$-4$&$-4$&$-10$&$2$&$2$\\\hhline{|=#=|=|=|=|=|=|=|=|}
$\mathfrak{p}$&$\mathfrak{p}_{37,3}$&$\mathfrak{p}_{37,4}$&$\mathfrak{p}_{61,1}$&$\mathfrak{p}_{61,2}$&$\mathfrak{p}_{61,3}$&$\mathfrak{p}_{73,1}$&$\mathfrak{p}_{73,3}$&$\mathfrak{p}_{73,4}$\\\hline
$a_\mathfrak{p}$&$2$&$8$&$2$&$-10$&$8$&$-16$&$14$&$-10$\\\hhline{|=#=|=|=|=|=|=|=|=|}
$\mathfrak{p}$&$\mathfrak{p}_{97,1}$&$\mathfrak{p}_{97,2}$&$\mathfrak{p}_{97,4}$&$\mathfrak{p}_{109,1}$&$\mathfrak{p}_{109,2}$&$\mathfrak{p}_{109,4}$&$\mathfrak{p}_{181,1}$&$\mathfrak{p}_{193,1}$\\\hline
$a_\mathfrak{p}$&$2$&$2$&$-4$&$2$&$2$&$-10$&$2$&$-22$\\\hhline{|=#=|=|=|=|=|=|=|=|}
$\mathfrak{p}$&$\mathfrak{p}_{229,3}$&$\mathfrak{p}_{17,1}$&$\mathfrak{p}_{313,1}$&$\mathfrak{p}_{313,4}$&$\mathfrak{p}_{373,1}$&$\mathfrak{p}_{409,1}$&$\mathfrak{p}_{1321,1}$&\\\hline
$a_\mathfrak{p}$&$-4$&$2$&$-10$&$14$&$-10$&$32$&$-10$&\\\hline
\end{tabular}\caption*{Eigenvalues $a_\mathfrak{p}$ of the Hecke operators $T_\mathfrak{p}$ on class 2041}\end{table}

\textbf{Class 2257}

\begin{itemize}
\item The residual representation attached to the corresponding elliptic curve has image isomorphic to $S_3$.
\item The primes $\{\mathfrak{p}_3,\mathfrak{p}_{13,1},\mathfrak{p}_{13,2},\mathfrak{p}_{13,4},\mathfrak{p}_{5,1}\}$ suffice to prove isomorphism of the residual representations.
\item The primes $\{\mathfrak{p}_{3}$, $\mathfrak{p}_{13,1}$, $\mathfrak{p}_{13,2}$, $\mathfrak{p}_{13,3}$, $\mathfrak{p}_{13,4}$, $\mathfrak{p}_{5,1}$, $\mathfrak{p}_{37,3}$, $\mathfrak{p}_{61,1}$, $\mathfrak{p}_{97,1}$, $\mathfrak{p}_{97,3}$, $\mathfrak{p}_{109,3}$, $\mathfrak{p}_{193,1}\}$ suffice to prove isomorphism of the full representations.
\end{itemize}

\begin{table}[h]\centering\scriptsize\begin{tabular}{|c||c|c|c|c|c|c|}\hline
$\mathfrak{p}$&$\mathfrak{p}_{3}$&$\mathfrak{p}_{13,1}$&$\mathfrak{p}_{13,2}$&$\mathfrak{p}_{13,3}$&$\mathfrak{p}_{13,4}$&$\mathfrak{p}_{5,1}$\\\hline
$a_\mathfrak{p}$&$-4$&$-1$&$1$&$-6$&$-3$&$1$\\\hhline{|=#=|=|=|=|=|=|}
$\mathfrak{p}$&$\mathfrak{p}_{37,3}$&$\mathfrak{p}_{61,1}$&$\mathfrak{p}_{97,1}$&$\mathfrak{p}_{97,3}$&$\mathfrak{p}_{109,3}$&$\mathfrak{p}_{193,1}$\\\hline
$a_\mathfrak{p}$&$-3$&$-12$&$0$&$-10$&$8$&$-10$\\\hline
\end{tabular}\caption*{Eigenvalues $a_\mathfrak{p}$ of the Hecke operators $T_\mathfrak{p}$ on class 2257}\end{table}

\textbf{Class 2452}

\begin{itemize}
\item The residual representation attached to the corresponding elliptic curve has image isomorphic to $S_3$.
\item The primes $\{\mathfrak{p}_3,\mathfrak{p}_{13,1},\mathfrak{p}_{13,3},\mathfrak{p}_{5,2}\}$ suffice to prove isomorphism of the residual representations.
\item The primes $\{\mathfrak{p}_{3}$, $\mathfrak{p}_{13,1}$, $\mathfrak{p}_{13,3}$, $\mathfrak{p}_{5,2}$, $\mathfrak{p}_{37,2}$, $\mathfrak{p}_{37,3}$, $\mathfrak{p}_{7,1}$, $\mathfrak{p}_{7,2}$, $\mathfrak{p}_{61,3}\}$ suffice to prove isomorphism of the full representations.
\end{itemize}

\begin{table}[h]\centering\scriptsize\begin{tabular}{|c||c|c|c|c|c|}\hline
$\mathfrak{p}$&$\mathfrak{p}_{3}$&$\mathfrak{p}_{13,1}$&$\mathfrak{p}_{13,3}$&$\mathfrak{p}_{5,2}$&$\mathfrak{p}_{37,2}$\\\hline
$a_\mathfrak{p}$&$1$&$-4$&$-4$&$8$&$2$\\\hhline{|=#=|=|=|=|=|}
$\mathfrak{p}$&$\mathfrak{p}_{37,3}$&$\mathfrak{p}_{7,1}$&$\mathfrak{p}_{7,2}$&$\mathfrak{p}_{61,3}$&\\\hline
$a_\mathfrak{p}$&$11$&$-4$&$-4$&$-10$&\\\hline
\end{tabular}\caption*{Eigenvalues $a_\mathfrak{p}$ of the Hecke operators $T_\mathfrak{p}$ on class 2452}\end{table}

\textbf{Class 2500a}

\begin{itemize}
\item The residual representation attached to the corresponding elliptic curve has image isomorphic to $S_3$.
\item The primes $\{\mathfrak{p}_{13,3},\mathfrak{p}_{13,4},\mathfrak{p}_{37,1}\}$ suffice to prove isomorphism of the residual representations.
\item The primes $\{\mathfrak{p}_{13,1}$, $\mathfrak{p}_{13,3}$, $\mathfrak{p}_{5,2}$, $\mathfrak{p}_{37,1}$, $\mathfrak{p}_{37,2}$, $\mathfrak{p}_{37,3}$, $\mathfrak{p}_{61,3}\}$ suffice to prove isomorphism of the full representations.
\end{itemize}

\newpage\begin{table}[h]\centering\scriptsize\begin{tabular}{|c||c|c|c|c|}\hline
$\mathfrak{p}$&$\mathfrak{p}_{13,1}$&$\mathfrak{p}_{13,3}$&$\mathfrak{p}_{13,4}$&$\mathfrak{p}_{5,2}$\\\hline
$a_\mathfrak{p}$&$4$&$-1$&$-1$&$1$ \\\hhline{|=#=|=|=|=|}
$\mathfrak{p}$&$\mathfrak{p}_{37,1}$&$\mathfrak{p}_{37,2}$&$\mathfrak{p}_{37,3}$&$\mathfrak{p}_{61,3}$\\\hline
$a_\mathfrak{p}$&$-7$&$-2$&$-7$&$-8$ \\\hline
\end{tabular}\caption*{Eigenvalues $a_\mathfrak{p}$ of the Hecke operators $T_\mathfrak{p}$ on class 2500a}\end{table}

\textbf{Class 2977}

\begin{itemize} 
\item The residual representation attached to the corresponding elliptic curve has trivial image.
\item The primes $\{\mathfrak{p}_{13,1},\mathfrak{p}_{13,2},\mathfrak{p}_{5,2}\}$ suffice to prove isomorphism of the residual representations.
\item The primes 
$\{\mathfrak{p}_{3}$, $\mathfrak{p}_{13,1}$, $\mathfrak{p}_{13,2}$, $\mathfrak{p}_{13,4}$, $\mathfrak{p}_{5,1}$, $\mathfrak{p}_{5,2}$, $\mathfrak{p}_{37,1}$, $\mathfrak{p}_{37,2}$, $\mathfrak{p}_{37,3}$, $\mathfrak{p}_{37,4}$, $\mathfrak{p}_{7,2}$, $\mathfrak{p}_{61,1}$, $\mathfrak{p}_{61,2}$, $\mathfrak{p}_{73,1}$, $\mathfrak{p}_{73,2}$, $\mathfrak{p}_{73,3}$, $\mathfrak{p}_{97,1}$, $\mathfrak{p}_{97,2}$, $\mathfrak{p}_{109,2}$, $\mathfrak{p}_{109,3}$, $\mathfrak{p}_{157,2}$, $\mathfrak{p}_{157,4}$, $\mathfrak{p}_{229,1}$, $\mathfrak{p}_{229,2}$, $\mathfrak{p}_{241,3}$, $\mathfrak{p}_{17,1}$, $\mathfrak{p}_{313,1}$, $\mathfrak{p}_{19,1}$, $\mathfrak{p}_{397,3}$, $\mathfrak{p}_{409,2}$, $\mathfrak{p}_{409,3}\}$ suffice to prove isomorphism of the full representations.
\end{itemize}

\begin{table}[h]\centering\scriptsize\begin{tabular}{|c||c|c|c|c|c|c|c|c|}\hline
$\mathfrak{p}$&$\mathfrak{p}_{3}$&$\mathfrak{p}_{13,1}$&$\mathfrak{p}_{13,2}$&$\mathfrak{p}_{13,4}$&$\mathfrak{p}_{5,1}$&$\mathfrak{p}_{5,2}$&$\mathfrak{p}_{37,1}$&$\mathfrak{p}_{37,2}$\\\hline
$a_\mathfrak{p}$&$4$&$2$&$-4$&$-4$&$2$&$-10$&$2$&$-10$\\\hhline{|=#=|=|=|=|=|=|=|=|}
$\mathfrak{p}$&$\mathfrak{p}_{37,3}$&$\mathfrak{p}_{37,4}$&$\mathfrak{p}_{7,2}$&$\mathfrak{p}_{61,1}$&$\mathfrak{p}_{61,2}$&$\mathfrak{p}_{73,1}$&$\mathfrak{p}_{73,2}$&$\mathfrak{p}_{73,3}$\\\hline
$a_\mathfrak{p}$&$-10$&$2$&$8$&$-10$&$2$&$-4$&$14$&$2$\\\hhline{|=#=|=|=|=|=|=|=|=|}
$\mathfrak{p}$&$\mathfrak{p}_{97,1}$&$\mathfrak{p}_{97,2}$&$\mathfrak{p}_{109,2}$&$\mathfrak{p}_{109,3}$&$\mathfrak{p}_{157,2}$&$\mathfrak{p}_{157,4}$&$\mathfrak{p}_{229,1}$&$\mathfrak{p}_{229,2}$\\\hline 
$a_\mathfrak{p}$&$2$&$14$&$14$&$-10$&$2$&$-16$&$8$&$2$\\\hhline{|=#=|=|=|=|=|=|=|=|}
$\mathfrak{p}$&$\mathfrak{p}_{241,3}$&$\mathfrak{p}_{17,1}$&$\mathfrak{p}_{313,1}$&$\mathfrak{p}_{19,1}$&$\mathfrak{p}_{397,3}$&$\mathfrak{p}_{409,2}$&$\mathfrak{p}_{409,3}$&\\\hline
$a_\mathfrak{p}$&$14$&$8$&$-10$&$2$&$2$&$14$&$2$&\\\hline
\end{tabular}\caption*{Eigenvalues $a_\mathfrak{p}$ of the Hecke operators $T_\mathfrak{p}$ on class 2977}\end{table}

\textbf{Class 3328}

\begin{itemize} 
\item The residual representation attached to the corresponding elliptic curve has trivial image.
\item Using the methods outlined in \textbf{[Jon14]}, one can immediately deduce that the residual representations are isomorphic.
\item The primes 
$\{\mathfrak{p}_{3}$, $\mathfrak{p}_{13,1}$, $\mathfrak{p}_{13,2}$, $\mathfrak{p}_{13,3}$, $\mathfrak{p}_{5,1}$, $\mathfrak{p}_{5,2}$, $\mathfrak{p}_{37,1}$, $\mathfrak{p}_{37,2}$, $\mathfrak{p}_{37,4}$, $\mathfrak{p}_{61,1}$, $\mathfrak{p}_{61,2}$, $\mathfrak{p}_{73,1}$, $\mathfrak{p}_{73,4}$, $\mathfrak{p}_{97,1}$, $\mathfrak{p}_{17,1}\}$ suffice to prove isomorphism of the full representations.
\end{itemize}

\begin{table}[h]\centering\scriptsize\begin{tabular}{|c||c|c|c|c|c|c|c|c|}\hline
$\mathfrak{p}$&$\mathfrak{p}_{3}$&$\mathfrak{p}_{13,1}$&$\mathfrak{p}_{13,2}$&$\mathfrak{p}_{13,3}$&$\mathfrak{p}_{5,1}$&$\mathfrak{p}_{5,2}$&$\mathfrak{p}_{37,1}$&$\mathfrak{p}_{37,2}$\\\hline
$a_\mathfrak{p}$&$2$&$-2$&$-2$&$6$&$2$&$-6$&$-10$&$-2$\\\hhline{|=#=|=|=|=|=|=|=|=|}
$\mathfrak{p}$&$\mathfrak{p}_{37,4}$&$\mathfrak{p}_{61,1}$&$\mathfrak{p}_{61,2}$&$\mathfrak{p}_{73,1}$&$\mathfrak{p}_{73,4}$&$\mathfrak{p}_{97,1}$&$\mathfrak{p}_{17,1}$&\\\hline
$a_\mathfrak{p}$&$6$&$-2$&$-2$&$10$&$10$&$18$&$2$&\\\hline
\end{tabular}\caption*{Eigenvalues $a_\mathfrak{p}$ of the Hecke operators $T_\mathfrak{p}$ on class 3328}\end{table}

\textbf{Class 3721b}

\begin{itemize} 
\item The residual representation attached to the corresponding elliptic curve has trivial image.
\item The primes $\{\mathfrak{p}_{13,1},\mathfrak{p}_{13,3}\}$ suffice to prove isomorphism of the residual representations.
\item The primes 
$\{\mathfrak{p}_{3}$, $\mathfrak{p}_{13,1}$, $\mathfrak{p}_{13,2}$, $\mathfrak{p}_{13,3}$, $\mathfrak{p}_{13,4}$, $\mathfrak{p}_{5,1}$, $\mathfrak{p}_{37,1}$, $\mathfrak{p}_{37,3}$, 
$\mathfrak{p}_{7,1}$, $\mathfrak{p}_{7,2}$, $\mathfrak{p}_{61,3}$, $\mathfrak{p}_{61,4}$, $\mathfrak{p}_{73,1}$, $\mathfrak{p}_{73,2}$, $\mathfrak{p}_{97,1}$, $\mathfrak{p}_{97,2}$, 
$\mathfrak{p}_{109,1}$, $\mathfrak{p}_{109,3}$, $\mathfrak{p}_{157,1}$, $\mathfrak{p}_{157,3}$, $\mathfrak{p}_{181,1}$, $\mathfrak{p}_{181,2}$, $\mathfrak{p}_{181,3}$, $\mathfrak{p}_{181,4}$, 
$\mathfrak{p}_{193,1}$, $\mathfrak{p}_{17,1}$, $\mathfrak{p}_{337,1}$, $\mathfrak{p}_{19,1}$, $\mathfrak{p}_{409,1}$, $\mathfrak{p}_{409,3}$, $\mathfrak{p}_{601,3}\}$ suffice to prove isomorphism of the full representations.
\end{itemize}

\begin{table}[h]\centering\scriptsize\begin{tabular}{|c||c|c|c|c|c|c|c|c|}\hline
$\mathfrak{p}$&$\mathfrak{p}_{3}$&$\mathfrak{p}_{13,1}$&$\mathfrak{p}_{13,2}$&$\mathfrak{p}_{13,3}$&$\mathfrak{p}_{13,4}$&$\mathfrak{p}_{5,1}$&$\mathfrak{p}_{37,1}$&$\mathfrak{p}_{37,3}$\\\hline
$a_\mathfrak{p}$&$-2$&$-4$&$-4$&$2$&$2$&$8$&$-10$&$-2$\\\hhline{|=#=|=|=|=|=|=|=|=|}
$\mathfrak{p}$&$\mathfrak{p}_{7,1}$&$\mathfrak{p}_{7,2}$&$\mathfrak{p}_{61,3}$&$\mathfrak{p}_{61,4}$&$\mathfrak{p}_{73,1}$&$\mathfrak{p}_{73,2}$&$\mathfrak{p}_{97,1}$&$\mathfrak{p}_{97,2}$\\\hline
$a_\mathfrak{p}$&$2$&$2$&$2$&$2$&$2$&$2$&$14$&$14$\\\hhline{|=#=|=|=|=|=|=|=|=|}
$\mathfrak{p}$&$\mathfrak{p}_{109,1}$&$\mathfrak{p}_{109,3}$&$\mathfrak{p}_{157,1}$&$\mathfrak{p}_{157,3}$&$\mathfrak{p}_{181,1}$&$\mathfrak{p}_{181,2}$&$\mathfrak{p}_{181,3}$&$\mathfrak{p}_{181,4}$\\\hline 
$a_\mathfrak{p}$&$-4$&$-4$&$-10$&$-10$&$2$&$2$&$2$&$2$\\\hhline{|=#=|=|=|=|=|=|=|=|}
$\mathfrak{p}$&$\mathfrak{p}_{193,1}$&$\mathfrak{p}_{17,1}$&$\mathfrak{p}_{337,1}$&$\mathfrak{p}_{19,1}$&$\mathfrak{p}_{409,1}$&$\mathfrak{p}_{409,3}$&$\mathfrak{p}_{601,3}$&\\\hline
$a_\mathfrak{p}$&$14$&$32$&$2$&$20$&$-22$&$-22$&$20$&\\\hline
\end{tabular}\caption*{Eigenvalues $a_\mathfrak{p}$ of the Hecke operators $T_\mathfrak{p}$ on class 3721b}\end{table}

\textbf{Class 3844}

\begin{itemize}
\item The residual representation attached to the corresponding elliptic curve has image isomorphic to $S_3$.
\item The primes $\{\mathfrak{p}_{3},\mathfrak{p}_{13,1},\mathfrak{p}_{13,3},\mathfrak{p}_{73,2}\}$ suffice to prove isomorphism of the residual representations.
\item The primes $\{\mathfrak{p}_{3}$, $\mathfrak{p}_{13,1}$, $\mathfrak{p}_{13,2}$, $\mathfrak{p}_{5,1}$, $\mathfrak{p}_{37,2}$, $\mathfrak{p}_{61,1}$, $\mathfrak{p}_{73,2}\}$ suffice to prove isomorphism of the full representations.
\end{itemize}

\begin{table}[h]\centering\scriptsize\begin{tabular}{|c||c|c|c|c|}\hline
$\mathfrak{p}$&$\mathfrak{p}_{3}$&$\mathfrak{p}_{13,1}$&$\mathfrak{p}_{13,2}$&$\mathfrak{p}_{13,3}$\\\hline
$a_\mathfrak{p}$&$-5$&$-1$&$-6$&$-6$\\\hhline{|=#=|=|=|=|}
$\mathfrak{p}$&$\mathfrak{p}_{5,1}$&$\mathfrak{p}_{37,2}$&$\mathfrak{p}_{61,1}$&$\mathfrak{p}_{73,2}$\\\hline
$a_\mathfrak{p}$&$1$&$3$&$-13$&$4$\\\hline
\end{tabular}\caption*{Eigenvalues $a_\mathfrak{p}$ of the Hecke operators $T_\mathfrak{p}$ on class 3844}\end{table}

\textbf{Class 4033a}

\begin{itemize} 
\item The residual representation attached to the corresponding elliptic curve has trivial image.
\item The prime $\mathfrak{p}_{13,1}$ suffices to prove isomorphism of the residual representations.
\item The primes 
$\{\mathfrak{p}_{3}$, $\mathfrak{p}_{13,1}$, $\mathfrak{p}_{13,2}$, $\mathfrak{p}_{13,3}$, $\mathfrak{p}_{13,4}$, $\mathfrak{p}_{5,1}$, $\mathfrak{p}_{5,2}$, $\mathfrak{p}_{37,1}$, 
$\mathfrak{p}_{37,2}$, $\mathfrak{p}_{37,3}$, $\mathfrak{p}_{7,1}$, $\mathfrak{p}_{61,1}$, $\mathfrak{p}_{61,3}$, $\mathfrak{p}_{61,4}$, $\mathfrak{p}_{73,1}$, $\mathfrak{p}_{73,2}$, 
$\mathfrak{p}_{73,3}$, $\mathfrak{p}_{73,4}$, $\mathfrak{p}_{97,1}$, $\mathfrak{p}_{97,2}$, $\mathfrak{p}_{97,3}$, $\mathfrak{p}_{109,1}$, $\mathfrak{p}_{109,2}$, $\mathfrak{p}_{157,4}$, 
$\mathfrak{p}_{277,1}$, $\mathfrak{p}_{17,1}$, $\mathfrak{p}_{17,2}$, $\mathfrak{p}_{337,4}$, $\mathfrak{p}_{349,2}$, $\mathfrak{p}_{373,1}$, $\mathfrak{p}_{409,3}\}$ suffice to prove isomorphism of the full representations.
\end{itemize}

\begin{table}[h]\centering\scriptsize\begin{tabular}{|c||c|c|c|c|c|c|c|c|}\hline
$\mathfrak{p}$&$\mathfrak{p}_{3}$&$\mathfrak{p}_{13,1}$&$\mathfrak{p}_{13,2}$&$\mathfrak{p}_{13,3}$&$\mathfrak{p}_{13,4}$&$\mathfrak{p}_{5,1}$&$\mathfrak{p}_{5,2}$&$\mathfrak{p}_{37,1}$\\\hline 
$a_\mathfrak{p}$&$4$&$2$&$-4$&$2$&$2$&$2$&$-4$&$2$\\\hhline{|=#=|=|=|=|=|=|=|=|}
$\mathfrak{p}$&$\mathfrak{p}_{37,2}$&$\mathfrak{p}_{37,3}$&$\mathfrak{p}_{7,1}$&$\mathfrak{p}_{61,1}$&$\mathfrak{p}_{61,3}$&$\mathfrak{p}_{61,4}$&$\mathfrak{p}_{73,1}$&$\mathfrak{p}_{73,2}$\\\hline 
$a_\mathfrak{p}$&$2$&$2$&$-4$&$2$&$14$&$-10$&$2$&$2$\\\hhline{|=#=|=|=|=|=|=|=|=|}
$\mathfrak{p}$&$\mathfrak{p}_{73,3}$&$\mathfrak{p}_{73,4}$&$\mathfrak{p}_{97,1}$&$\mathfrak{p}_{97,2}$&$\mathfrak{p}_{97,3}$&$\mathfrak{p}_{109,1}$&$\mathfrak{p}_{109,2}$&$\mathfrak{p}_{157,4}$\\\hline 
$a_\mathfrak{p}$&$2$&$2$&$-16$&$2$&$14$&$2$&$-16$&$-4$\\\hhline{|=#=|=|=|=|=|=|=|=|}
$\mathfrak{p}$&$\mathfrak{p}_{277,1}$&$\mathfrak{p}_{17,1}$&$\mathfrak{p}_{17,2}$&$\mathfrak{p}_{337,4}$&$\mathfrak{p}_{349,2}$&$\mathfrak{p}_{373,1}$&$\mathfrak{p}_{409,3}$&\\\hline
$a_\mathfrak{p}$&$-10$&$2$&$2$&$2$&$32$&$14$&$-10$&\\\hline
\end{tabular}\caption*{Eigenvalues $a_\mathfrak{p}$ of the Hecke operators $T_\mathfrak{p}$ on class 4033a}\end{table}

\newpage\textbf{Class 4033b}

\begin{itemize} 
\item The residual representation attached to the corresponding elliptic curve has trivial image.
\item The primes $\{\mathfrak{p}_3,\mathfrak{p}_{13,2}\}$ suffice to prove isomorphism of the residual representations.
\item The primes 
$\{\mathfrak{p}_{3}$, $\mathfrak{p}_{13,1}$, $\mathfrak{p}_{13,2}$, $\mathfrak{p}_{13,3}$, $\mathfrak{p}_{13,4}$, $\mathfrak{p}_{5,1}$, $\mathfrak{p}_{5,2}$, $\mathfrak{p}_{37,1}$, 
$\mathfrak{p}_{37,2}$, $\mathfrak{p}_{37,3}$, $\mathfrak{p}_{7,1}$, $\mathfrak{p}_{61,1}$, $\mathfrak{p}_{61,2}$, $\mathfrak{p}_{61,3}$, $\mathfrak{p}_{61,4}$, $\mathfrak{p}_{73,1}$, 
$\mathfrak{p}_{73,2}$, $\mathfrak{p}_{73,3}$, $\mathfrak{p}_{73,4}$, $\mathfrak{p}_{97,1}$, $\mathfrak{p}_{97,2}$, $\mathfrak{p}_{97,3}$, $\mathfrak{p}_{109,1}$, $\mathfrak{p}_{109,2}$, 
$\mathfrak{p}_{109,4}$, $\mathfrak{p}_{157,4}$, $\mathfrak{p}_{17,2}$, $\mathfrak{p}_{313,4}$, $\mathfrak{p}_{373,1}$, $\mathfrak{p}_{409,2}$, $\mathfrak{p}_{409,3}\}$ suffice to prove isomorphism of the full representations.
\end{itemize}

\begin{table}[h]\centering\scriptsize\begin{tabular}{|c||c|c|c|c|c|c|c|c|}\hline
$\mathfrak{p}$&$\mathfrak{p}_{3}$&$\mathfrak{p}_{13,1}$&$\mathfrak{p}_{13,2}$&$\mathfrak{p}_{13,3}$&$\mathfrak{p}_{13,4}$&$\mathfrak{p}_{5,1}$&$\mathfrak{p}_{5,2}$&$\mathfrak{p}_{37,1}$\\\hline 
$a_\mathfrak{p}$&$-2$&$-4$&$2$&$2$&$2$&$-4$&$2$&$2$\\\hhline{|=#=|=|=|=|=|=|=|=|}
$\mathfrak{p}$&$\mathfrak{p}_{37,2}$&$\mathfrak{p}_{37,3}$&$\mathfrak{p}_{7,1}$&$\mathfrak{p}_{61,1}$&$\mathfrak{p}_{61,2}$&$\mathfrak{p}_{61,3}$&$\mathfrak{p}_{61,4}$&$\mathfrak{p}_{73,1}$\\\hline 
$a_\mathfrak{p}$&$2$&$2$&$-10$&$8$&$-4$&$8$&$2$&$2$\\\hhline{|=#=|=|=|=|=|=|=|=|}
$\mathfrak{p}$&$\mathfrak{p}_{73,2}$&$\mathfrak{p}_{73,3}$&$\mathfrak{p}_{73,4}$&$\mathfrak{p}_{97,1}$&$\mathfrak{p}_{97,2}$&$\mathfrak{p}_{97,3}$&$\mathfrak{p}_{109,1}$&$\mathfrak{p}_{109,2}$\\\hline 
$a_\mathfrak{p}$&$2$&$2$&$-16$&$-10$&$8$&$-10$&$2$&$2$\\\hhline{|=#=|=|=|=|=|=|=|=|}
$\mathfrak{p}$&$\mathfrak{p}_{109,4}$&$\mathfrak{p}_{157,4}$&$\mathfrak{p}_{17,2}$&$\mathfrak{p}_{313,4}$&$\mathfrak{p}_{373,1}$&$\mathfrak{p}_{409,2}$&$\mathfrak{p}_{409,3}$&\\\hline
$a_\mathfrak{p}$&$2$&$-10$&$20$&$26$&$-34$&$2$&$14$&\\\hline
\end{tabular}\caption*{Eigenvalues $a_\mathfrak{p}$ of the Hecke operators $T_\mathfrak{p}$ on class 4033b}\end{table}

\textbf{Class 4057}

\begin{itemize}
\item The residual representation attached to the corresponding elliptic curve has image isomorphic to $S_3$.
\item The primes $\{\mathfrak{p}_{13,1},\mathfrak{p}_{13,3},\mathfrak{p}_{5,2}\}$ suffice to prove isomorphism of the residual representations.
\item The primes $\{\mathfrak{p}_{3}$, $\mathfrak{p}_{13,1}$, $\mathfrak{p}_{13,2}$, $\mathfrak{p}_{13,3}$, $\mathfrak{p}_{5,1}$, $\mathfrak{p}_{5,2}$, $\mathfrak{p}_{37,2}$, $\mathfrak{p}_{61,2}$, $\mathfrak{p}_{61,4}
\}$ suffice to prove isomorphism of the full representations.
\end{itemize}

\begin{table}[h]\centering\scriptsize\begin{tabular}{|c||c|c|c|c|c|}\hline
$\mathfrak{p}$&$\mathfrak{p}_{3}$&$\mathfrak{p}_{13,1}$&$\mathfrak{p}_{13,2}$&$\mathfrak{p}_{13,3}$&$\mathfrak{p}_{5,1}$\\\hline
$a_\mathfrak{p}$&$-2$&$-4$&$-1$&$-4$&$-5$\\\hhline{|=#=|=|=|=|=|}
$\mathfrak{p}$&$\mathfrak{p}_{5,2}$&$\mathfrak{p}_{37,2}$&$\mathfrak{p}_{61,2}$&$\mathfrak{p}_{61,4}$&\\\hline
$a_\mathfrak{p}$&$-2$&$4$&$-13$&$10$&\\\hline
\end{tabular}\caption*{Eigenvalues $a_\mathfrak{p}$ of the Hecke operators $T_\mathfrak{p}$ on class 4057}\end{table}

\textbf{Class 4069}

\begin{itemize}
\item The residual representation attached to the corresponding elliptic curve has image isomorphic to $S_3$.
\item The primes $\{\mathfrak{p}_{13,1},\mathfrak{p}_{13,2},\mathfrak{p}_{13,4},\mathfrak{p}_{5,2},\mathfrak{p}_{37,2}\}$ suffice prove isomorphism of the residual representations.
\item The primes $\{\mathfrak{p}_{13,1}$, $\mathfrak{p}_{13,2}$, $\mathfrak{p}_{13,4}$, $\mathfrak{p}_{5,2}$, $\mathfrak{p}_{37,1}$, $\mathfrak{p}_{37,2}$, $\mathfrak{p}_{37,3}$, $\mathfrak{p}_{61,3}$, $\mathfrak{p}_{73,2}$, $\mathfrak{p}_{97,2}$, $\mathfrak{p}_{97,3}$, $\mathfrak{p}_{109,3}\}$ suffice to prove isomorphism of the full representations.
\end{itemize}

\begin{table}[h]\centering\scriptsize\begin{tabular}{|c||c|c|c|c|c|c|}\hline
$\mathfrak{p}$&$\mathfrak{p}_{13,1}$&$\mathfrak{p}_{13,2}$&$\mathfrak{p}_{13,4}$&$\mathfrak{p}_{5,2}$&$\mathfrak{p}_{37,1}$&$\mathfrak{p}_{37,2}$\\\hline
$a_\mathfrak{p}$&$-3$&$1$&$-5$&$1$&$-7$&$-10$ \\\hhline{|=#=|=|=|=|=|=|}
$\mathfrak{p}$&$\mathfrak{p}_{37,3}$&$\mathfrak{p}_{61,3}$&$\mathfrak{p}_{73,2}$&$\mathfrak{p}_{97,2}$&$\mathfrak{p}_{97,3}$&$\mathfrak{p}_{109,3}$\\\hline
$\mathfrak{p}$&$-2$&$-12$&$14$&$10$&$8$&$-10$ \\\hline
\end{tabular}\caption*{Eigenvalues $a_\mathfrak{p}$ of the Hecke operators $T_\mathfrak{p}$ on class 4069}\end{table}

\newpage\textbf{Class 4225b}

\begin{itemize} 
\item The residual representation attached to the corresponding elliptic curve has trivial image.
\item The primes $\{\mathfrak{p}_3,\mathfrak{p}_{13,3}\}$ suffice to prove isomorphism of the residual representations.
\item The primes 
$\{\mathfrak{p}_{3}$, $\mathfrak{p}_{13,1}$, $\mathfrak{p}_{13,3}$, $\mathfrak{p}_{13,4}$, $\mathfrak{p}_{5,2}$, $\mathfrak{p}_{37,1}$, $\mathfrak{p}_{37,2}$, $\mathfrak{p}_{37,3}$, 
$\mathfrak{p}_{37,4}$, $\mathfrak{p}_{7,1}$, $\mathfrak{p}_{7,2}$, $\mathfrak{p}_{61,1}$, $\mathfrak{p}_{61,3}$, $\mathfrak{p}_{61,4}$, $\mathfrak{p}_{73,1}$, $\mathfrak{p}_{73,2}$, 
$\mathfrak{p}_{73,3}$, $\mathfrak{p}_{97,2}$, $\mathfrak{p}_{97,3}$, $\mathfrak{p}_{109,3}$, $\mathfrak{p}_{11,1}$, $\mathfrak{p}_{157,2}$, $\mathfrak{p}_{157,4}$, $\mathfrak{p}_{181,4}$, 
$\mathfrak{p}_{193,4}$, $\mathfrak{p}_{229,1}$, $\mathfrak{p}_{229,2}$, $\mathfrak{p}_{17,1}$, $\mathfrak{p}_{313,3}$, $\mathfrak{p}_{409,1}$, $\mathfrak{p}_{409,2}\}$ suffice to prove isomorphism of the full representations.
\end{itemize}

\begin{table}[h]\centering\scriptsize\begin{tabular}{|c||c|c|c|c|c|c|c|c|}\hline
$\mathfrak{p}$&$\mathfrak{p}_{3}$&$\mathfrak{p}_{13,1}$&$\mathfrak{p}_{13,3}$&$\mathfrak{p}_{13,4}$&$\mathfrak{p}_{5,2}$&$\mathfrak{p}_{37,1}$&$\mathfrak{p}_{37,2}$&$\mathfrak{p}_{37,3}$\\\hline 
$a_\mathfrak{p}$&$-2$&$-4$&$-2$&$-6$&$4$&$0$&$-2$&$-6$\\\hhline{|=#=|=|=|=|=|=|=|=|}
$\mathfrak{p}$&$\mathfrak{p}_{37,4}$&$\mathfrak{p}_{7,1}$&$\mathfrak{p}_{7,2}$&$\mathfrak{p}_{61,1}$&$\mathfrak{p}_{61,3}$&$\mathfrak{p}_{61,4}$&$\mathfrak{p}_{73,1}$&$\mathfrak{p}_{73,2}$\\\hline 
$a_\mathfrak{p}$&$-2$&$-6$&$2$&$-6$&$-10$&$-2$&$8$&$2$ \\\hhline{|=#=|=|=|=|=|=|=|=|}
$\mathfrak{p}$&$\mathfrak{p}_{73,3}$&$\mathfrak{p}_{97,2}$&$\mathfrak{p}_{97,3}$&$\mathfrak{p}_{109,3}$&$\mathfrak{p}_{11,1}$&$\mathfrak{p}_{157,2}$&$\mathfrak{p}_{157,4}$&$\mathfrak{p}_{181,4}$\\\hline 
$a_\mathfrak{p}$&$-14$&$-2$&$6$&$14$&$-12$&$-6$&$4$&$-10$ \\\hhline{|=#=|=|=|=|=|=|=|=|}
$\mathfrak{p}$&$\mathfrak{p}_{193,4}$&$\mathfrak{p}_{229,1}$&$\mathfrak{p}_{229,2}$&$\mathfrak{p}_{17,1}$&$\mathfrak{p}_{313,3}$&$\mathfrak{p}_{409,1}$&$\mathfrak{p}_{409,2}$&\\\hline
$a_\mathfrak{p}$&$-2$&$-14$&$20$&$-24$&$-10$&$10$&$-38$&\\\hline
\end{tabular}\caption*{Eigenvalues $a_\mathfrak{p}$ of the Hecke operators $T_\mathfrak{p}$ on class 4225b}\end{table}

\textbf{Class 4516}

\begin{itemize}
\item The residual representation attached to the corresponding elliptic curve has image isomorphic to $S_3$.
\item The primes $\{\mathfrak{p}_{3},\mathfrak{p}_{13,2},\mathfrak{p}_{13,3},\mathfrak{p}_{37,3}\}$ suffice to prove isomorphism of the residual representations.
\item The primes $\{\mathfrak{p}_{3}$, $\mathfrak{p}_{13,1}$, $\mathfrak{p}_{13,2}$, $\mathfrak{p}_{13,3}$, $\mathfrak{p}_{5,1}$, $\mathfrak{p}_{5,2}$, $\mathfrak{p}_{37,1}$, $\mathfrak{p}_{37,3}$, $\mathfrak{p}_{7,1}$, $\mathfrak{p}_{73,1}\}$ suffice to prove isomorphism of the full representations.
\end{itemize}

\begin{table}[h]\centering\scriptsize\begin{tabular}{|c||c|c|c|c|c|}\hline
$\mathfrak{p}$&$\mathfrak{p}_{3}$&$\mathfrak{p}_{13,1}$&$\mathfrak{p}_{13,2}$&$\mathfrak{p}_{13,3}$&$\mathfrak{p}_{5,1}$\\\hline
$a_\mathfrak{p}$&$5$&$4$&$-1$&$-6$&$-4$ \\\hhline{|=#=|=|=|=|=|}
$\mathfrak{p}$&$\mathfrak{p}_{5,2}$&$\mathfrak{p}_{37,1}$&$\mathfrak{p}_{37,3}$&$\mathfrak{p}_{7,1}$&$\mathfrak{p}_{73,1}$\\\hline
$a_\mathfrak{p}$&$6$&$-12$&$3$&$-10$&$-11$\\\hline
\end{tabular}\caption*{Eigenvalues $a_\mathfrak{p}$ of the Hecke operators $T_\mathfrak{p}$ on class 4516}\end{table}

\newpage\textbf{Class 4672}

\begin{itemize} 
\item The residual representation attached to the corresponding elliptic curve has trivial image.
\item Using the methods outlined in \textbf{[Jon14]}, one can immediately deduce that the residual representations are isomorphic.
\item The primes 
$\{\mathfrak{p}_{3}$, $\mathfrak{p}_{13,1}$, $\mathfrak{p}_{13,2}$, $\mathfrak{p}_{13,3}$, $\mathfrak{p}_{13,4}$, $\mathfrak{p}_{5,1}$, $\mathfrak{p}_{5,2}$, $\mathfrak{p}_{37,3}$, 
$\mathfrak{p}_{61,3}$, $\mathfrak{p}_{61,4}$, $\mathfrak{p}_{73,2}$, $\mathfrak{p}_{73,3}$, $\mathfrak{p}_{97,2}$, $\mathfrak{p}_{109,1}$, $\mathfrak{p}_{17,2}\}$ suffice to prove isomorphism of the full representations.
\end{itemize}

\begin{table}[h]\centering\scriptsize\begin{tabular}{|c||c|c|c|c|c|c|c|c|}\hline
$\mathfrak{p}$&$\mathfrak{p}_{3}$&$\mathfrak{p}_{13,1}$&$\mathfrak{p}_{13,2}$&$\mathfrak{p}_{13,3}$&$\mathfrak{p}_{13,4}$&$\mathfrak{p}_{5,1}$&$\mathfrak{p}_{5,2}$&$\mathfrak{p}_{37,3}$\\\hline 
$a_\mathfrak{p}$&$2$&$-2$&$-2$&$-2$&$6$&$-6$&$2$&$6$\\\hhline{|=#=|=|=|=|=|=|=|=|}
$\mathfrak{p}$&$\mathfrak{p}_{61,3}$&$\mathfrak{p}_{61,4}$&$\mathfrak{p}_{73,2}$&$\mathfrak{p}_{73,3}$&$\mathfrak{p}_{97,2}$&$\mathfrak{p}_{109,1}$&$\mathfrak{p}_{17,2}$&\\\hline
$a_\mathfrak{p}$&$-2$&$-10$&$-6$&$10$&$2$&$14$&$-30$&\\\hline
\end{tabular}\caption*{Eigenvalues $a_\mathfrak{p}$ of the Hecke operators $T_\mathfrak{p}$ on class 4672}\end{table}

\textbf{Class 4852}

\begin{itemize}
\item The residual representation attached to the corresponding elliptic curve has image isomorphic to $S_3$.
\item The primes $\{\mathfrak{p}_{3},\mathfrak{p}_{13,1},\mathfrak{p}_{13,2},\mathfrak{p}_{13,3}\}$ suffice to prove isomorphism of the residual representations.
\item The primes $\{\mathfrak{p}_{3}$, $\mathfrak{p}_{13,1}$, $\mathfrak{p}_{13,2}$, $\mathfrak{p}_{5,1}$, $\mathfrak{p}_{5,2}$, $\mathfrak{p}_{37,3}$, $\mathfrak{p}_{61,2}$, $\mathfrak{p}_{61,3}$, $\mathfrak{p}_{73,2}$, $\mathfrak{p}_{97,3}\}$ suffice to prove isomorphism of the full representations.
\end{itemize}

\begin{table}[h]\centering\scriptsize\begin{tabular}{|c||c|c|c|c|c|c|}\hline
$\mathfrak{p}$&$\mathfrak{p}_{3}$&$\mathfrak{p}_{13,1}$&$\mathfrak{p}_{13,2}$&$\mathfrak{p}_{13,3}$&$\mathfrak{p}_{5,1}$&$\mathfrak{p}_{5,2}$\\\hline
$a_\mathfrak{p}$&$-3$&$-1$&$-7$&$-2$&$3$&$-8$\\\hhline{|=#=|=|=|=|=|=|}
$\mathfrak{p}$&$\mathfrak{p}_{37,3}$&$\mathfrak{p}_{61,2}$&$\mathfrak{p}_{61,3}$&$\mathfrak{p}_{73,2}$&$\mathfrak{p}_{97,3}$&\\\hline
$a_\mathfrak{p}$&$2$&$-12$&$4$&$14$&$6$&\\\hline
\end{tabular}\caption*{Eigenvalues $a_\mathfrak{p}$ of the Hecke operators $T_\mathfrak{p}$ on class 4852}\end{table}

\textbf{Class 5317}

\begin{itemize} 
\item The residual representation attached to the corresponding elliptic curve has trivial image.
\item The primes $\{\mathfrak{p}_3,\mathfrak{p}_{13,2}\}$ suffice to prove isomorphism of the residual representations.
\item The primes 
$\{\mathfrak{p}_{3}$, $\mathfrak{p}_{13,1}$, $\mathfrak{p}_{13,2}$, $\mathfrak{p}_{13,3}$, $\mathfrak{p}_{5,1}$, $\mathfrak{p}_{5,2}$, $\mathfrak{p}_{37,1}$, $\mathfrak{p}_{37,2}$, 
$\mathfrak{p}_{37,4}$, $\mathfrak{p}_{7,2}$, $\mathfrak{p}_{61,1}$, $\mathfrak{p}_{61,2}$, $\mathfrak{p}_{61,3}$, $\mathfrak{p}_{73,1}$, $\mathfrak{p}_{73,2}$, $\mathfrak{p}_{73,4}$, 
$\mathfrak{p}_{97,1}$, $\mathfrak{p}_{97,3}$, $\mathfrak{p}_{97,4}$, $\mathfrak{p}_{109,3}$, $\mathfrak{p}_{157,2}$, $\mathfrak{p}_{157,3}$, $\mathfrak{p}_{181,1}$, $\mathfrak{p}_{181,2}$, 
$\mathfrak{p}_{193,2}$, $\mathfrak{p}_{277,4}$, $\mathfrak{p}_{17,1}$, $\mathfrak{p}_{313,2}$, $\mathfrak{p}_{373,3}$, $\mathfrak{p}_{409,1}$, $\mathfrak{p}_{457,1}\}$ suffice to prove isomorphism of the full representations.
\end{itemize}

\newpage\begin{table}[h]\centering\scriptsize\begin{tabular}{|c||c|c|c|c|c|c|c|c|}\hline
$\mathfrak{p}$&$\mathfrak{p}_{3}$&$\mathfrak{p}_{13,1}$&$\mathfrak{p}_{13,2}$&$\mathfrak{p}_{13,3}$&$\mathfrak{p}_{5,1}$&$\mathfrak{p}_{5,2}$&$\mathfrak{p}_{37,1}$&$\mathfrak{p}_{37,2}$\\\hline
$a_\mathfrak{p}$&$2$&$-2$&$6$&$-2$&$2$&$2$&$6$&$-10$\\\hhline{|=#=|=|=|=|=|=|=|=|}
$\mathfrak{p}$&$\mathfrak{p}_{37,4}$&$\mathfrak{p}_{7,2}$&$\mathfrak{p}_{61,1}$&$\mathfrak{p}_{61,2}$&$\mathfrak{p}_{61,3}$&$\mathfrak{p}_{73,1}$&$\mathfrak{p}_{73,2}$&$\mathfrak{p}_{73,4}$\\\hline
$a_\mathfrak{p}$&$6$&$2$&$-2$&$-10$&$-10$&$10$&$10$&$10$\\\hhline{|=#=|=|=|=|=|=|=|=|}
$\mathfrak{p}$&$\mathfrak{p}_{97,1}$&$\mathfrak{p}_{97,3}$&$\mathfrak{p}_{97,4}$&$\mathfrak{p}_{109,3}$&$\mathfrak{p}_{157,2}$&$\mathfrak{p}_{157,3}$&$\mathfrak{p}_{181,1}$&$\mathfrak{p}_{181,2}$\\\hline
$a_\mathfrak{p}$&$-14$&$2$&$-6$&$-2$&$-2$&$-2$&$-2$&$22$\\\hhline{|=#=|=|=|=|=|=|=|=|}
$\mathfrak{p}$&$\mathfrak{p}_{193,2}$&$\mathfrak{p}_{277,4}$&$\mathfrak{p}_{17,1}$&$\mathfrak{p}_{313,2}$&$\mathfrak{p}_{373,3}$&$\mathfrak{p}_{409,1}$&$\mathfrak{p}_{457,1}$&\\\hline
$a_\mathfrak{p}$&$2$&$6$&$34$&$-22$&$22$&$10$&$-22$&\\\hline
\end{tabular}\caption*{Eigenvalues $a_\mathfrak{p}$ of the Hecke operators $T_\mathfrak{p}$ on class 5317}\end{table}

\textbf{Class 5473}

\begin{itemize} 
\item The residual representation attached to the corresponding elliptic curve has trivial image.
\item The primes $\{\mathfrak{p}_3,\mathfrak{p}_{13,1}\}$ suffice to prove isomorphism of the residual representations.
\item The primes 
$\{\mathfrak{p}_{3}$, $\mathfrak{p}_{13,1}$, $\mathfrak{p}_{13,2}$, $\mathfrak{p}_{13,3}$, $\mathfrak{p}_{5,1}$, $\mathfrak{p}_{5,2}$, $\mathfrak{p}_{37,1}$, $\mathfrak{p}_{37,2}$, 
$\mathfrak{p}_{37,3}$, $\mathfrak{p}_{37,4}$, $\mathfrak{p}_{61,1}$, $\mathfrak{p}_{61,2}$, $\mathfrak{p}_{73,1}$, $\mathfrak{p}_{73,2}$, $\mathfrak{p}_{73,4}$, $\mathfrak{p}_{97,1}$, 
$\mathfrak{p}_{97,3}$, $\mathfrak{p}_{97,4}$, $\mathfrak{p}_{109,1}$, $\mathfrak{p}_{109,3}$, $\mathfrak{p}_{109,4}$, $\mathfrak{p}_{157,1}$, $\mathfrak{p}_{157,2}$, $\mathfrak{p}_{181,4}$, 
$\mathfrak{p}_{193,2}$, $\mathfrak{p}_{17,1}$, $\mathfrak{p}_{313,1}$, $\mathfrak{p}_{313,3}$, $\mathfrak{p}_{349,1}$, $\mathfrak{p}_{409,2}$, $\mathfrak{p}_{457,4}\}$ suffice to prove isomorphism of the full representations.
\end{itemize}

\begin{table}[h]\centering\scriptsize\begin{tabular}{|c||c|c|c|c|c|c|c|c|}\hline
$\mathfrak{p}$&$\mathfrak{p}_{3}$&$\mathfrak{p}_{13,1}$&$\mathfrak{p}_{13,2}$&$\mathfrak{p}_{13,3}$&$\mathfrak{p}_{5,1}$&$\mathfrak{p}_{5,2}$&$\mathfrak{p}_{37,1}$&$\mathfrak{p}_{37,2}$\\\hline 
$a_\mathfrak{p}$&$-2$&$2$&$2$&$4$&$2$&$8$&$2$&$-4$\\\hhline{|=#=|=|=|=|=|=|=|=|}
$\mathfrak{p}$&$\mathfrak{p}_{37,3}$&$\mathfrak{p}_{37,4}$&$\mathfrak{p}_{61,1}$&$\mathfrak{p}_{61,2}$&$\mathfrak{p}_{73,1}$&$\mathfrak{p}_{73,2}$&$\mathfrak{p}_{73,4}$&$\mathfrak{p}_{97,1}$\\\hline 
$a_\mathfrak{p}$&$2$&$2$&$14$&$8$&$14$&$14$&$14$&$8$ \\\hhline{|=#=|=|=|=|=|=|=|=|}
$\mathfrak{p}$&$\mathfrak{p}_{97,3}$&$\mathfrak{p}_{97,4}$&$\mathfrak{p}_{109,1}$&$\mathfrak{p}_{109,3}$&$\mathfrak{p}_{109,4}$&$\mathfrak{p}_{157,1}$&$\mathfrak{p}_{157,2}$&$\mathfrak{p}_{181,4}$\\\hline 
$a_\mathfrak{p}$&$8$&$2$&$2$&$14$&$-10$&$-4$&$-22$&$2$ \\\hhline{|=#=|=|=|=|=|=|=|=|}
$\mathfrak{p}$&$\mathfrak{p}_{193,2}$&$\mathfrak{p}_{17,1}$&$\mathfrak{p}_{313,1}$&$\mathfrak{p}_{313,3}$&$\mathfrak{p}_{349,1}$&$\mathfrak{p}_{409,2}$&$\mathfrak{p}_{457,4}$&\\\hline 
$a_\mathfrak{p}$&$-22$&$-10$&$8$&$-10$&$2$&$-22$&$-12$&\\\hline
\end{tabular}\caption*{Eigenvalues $a_\mathfrak{p}$ of the Hecke operators $T_\mathfrak{p}$ on class 5473}\end{table}

\section*{Acknowledgements}

I would like to thank Tobias Berger for suggesting this topic, and for his continued support and encouragement throughout, and to extend my gratitude to Paul Gunnells and Dan Yasaki for their helpful advice. Finally, I would like to thank the EPSRC, whose funding made this study possible.

\newpage

\section*{Appendix A. Proving Isomorphism of Residual Representations}

We present an algorithm (originally discussed in \textbf{[DGP10]}) which, given rational $2$-adic representations $\rho_E$ and $\rho_\pi$ attached to an elliptic curve $E$ over a CM field $F$ and a cuspidal automorphic representation $\pi$ of $\mathrm{Res}_{F/\mq}(\mathrm{GL}_2)$ respectively, allows one to determine isomorphism of the corresponding residual representations. This is done by considering the corresponding fields \begin{equation*}L_E := \overline{F}^{\mathrm{ker}(\widetilde{\rho_E})}~\mathrm{and}~L_\pi := \overline{F}^{\mathrm{ker}(\widetilde{\rho_\pi})}.\end{equation*}

If $L_E$ and $L_\pi$ are isomorphic over $F$, then the residual representations $\widetilde{\rho_E}$ and $\widetilde{\rho_\pi}$ are isomorphic. Indeed, these representations factor through the Galois groups $\mathrm{Gal}(L_E/F)$ and $\mathrm{Gal}(L_\pi/F)$ respectively, and an isomorphism between the fields implies that $\widetilde{\rho_E}(\mathrm{Frob}_\mathfrak{p})$ and $\widetilde{\rho_\pi}(\mathrm{Frob}_\mathfrak{p})$ must have the same order for all primes $\mathfrak{p}$ of $F$. Since the order of an element in $\mathrm{GL}_2(\mf_2)$ determines its characteristic polynomial, a theorem of Brauer and Nesbitt then implies isomorphism of the residual representations.
 
Given an elliptic curve $E$, it is a straightforward task to determine the image of $\widetilde{\rho_E}$, as it is defined by the $2$-torsion of $E$, which is simple to compute. Thus we shall assume knowledge of $L_E$, and reduce our task to determining whether $L_\pi$ is isomorphic to it. Now, $L_E$ and $L_\pi$ are unramified at all places $v$ of $F$ for which the representations $\widetilde{\rho_E}$ and $\widetilde{\rho_\pi}$ are unramified respectively, so we can assume that both fields are unramified at all places not dividing $2\mathfrak{n}$.

Since $\mathrm{Gal}(L_\pi/F) \simeq \mathrm{Im}(\widetilde{\rho_\pi})$  is a subgroup of $S_3$, we can build the field $L_\pi$ up from $F$ by extensions of degrees $2$ and $3$. Since these extensions are necessarily abelian, we are lead to consider class field theory.

Given an extension $K$ of $F$, and a list $p_1,\ldots,p_n$ of primes, let $r = p_1\ldots p_n$, and define a \it{modulus} \begin{equation*}\mathfrak{m}_{K,r} := \underset{\mathfrak{p}|2\mathfrak{n}}{\displaystyle\prod}~\mathfrak{p}^{e(\mathfrak{p})},\end{equation*} where \begin{equation*}e(\mathfrak{p}) = \left\{\begin{array}{cc}\lfloor \frac{p_ie(\mathfrak{p}/p_i)}{p_i-1}\rfloor + 1;~&\mathrm{if}~\mathfrak{p}|p_i,\\1;~&\mathrm{if~}\mathfrak{p}\nmid r.\end{array}\right.\end{equation*}

It can be shown (see, for example, \textbf{[Coh00], Propositions 3.3.21} and \textbf{3.3.22}) that if $L$ is a cyclic extension of such a field $K$, of order $p$, then $L$ corresponds to a character of order $p$ of the ray class group $Cl(\mathcal{O}_K,\mathfrak{m}_{K,r})$ for any $r$ which is divisible by $p$. We choose to work with \it{additive} characters, so that the group of order $p$ characters forms a $\mz/p\mz$-vector space. Fixing an isomorphism $\mathrm{Gal}(L/K) \simeq \mz/p\mz$, we can define the corresponding character $\chi_L$ on a prime ideal $\mathfrak{p}$ of $K$ not in $S$ by choosing a Frobenius element $\mathrm{Frob}_{L/K}(\mathfrak{p})$, and defining $\chi_L(\mathfrak{p})$ to be the image of this element in $\mz/p\mz$. We can then extend this to the entire ray class group in the obvious manner.

\newpage Observe that, given a residual representation $\widetilde{\rho}$ whose image lies in $\mathrm{GL}_2(\mf_2)$, with fixed field $\Lrho := \overline{F}^{\mathrm{ker}(\widetilde{\rho})}$, we can construct the following tower of field extensions: 

\begin{align*}
\xymatrix{
\Lrho \ar@{-}[d]  \ar@{-}@/_1pc/[dd]_{\mathrm{Gal}(\Lrho/F)} \ar@{-}@/^1pc/[d]^{\mathrm{Gal}(\Lrho/K)}\\
K \ar@{-}[d] \ar@{-}@/^1pc/[d]^{\mathrm{Gal}(K/F)}\\
F
}
\end{align*}
where \it{either} $\Lrho/F$ is trivial, \it{or} there exists an intermediate field $K$ such that $\Lrho/K$ is a cubic Galois extension, in which case $K/F$ is an extension of degree at most $2$.

Since $\widetilde{\rho}$ factors through $\mathrm{Gal}(\Lrho/F)$, we can deduce information about the image of $\widetilde{\rho}(\mathrm{Frob}_\mathfrak{p})$ for a prime $\mathfrak{p}$ in $F$ from the splitting behaviour of $\mathfrak{p}$. In particular:
\begin{itemize}
\item If $K$ is a quadratic extension of $F$, and $\mathfrak{p}$ is inert in $K$, then $\mathrm{Frob}_\mathfrak{p}|_{\Lrho}$ has order $2$ in $\mathrm{Gal}(\Lrho/F)$, and thus $\widetilde{\rho}(\mathrm{Frob}_\mathfrak{p})$ has even trace.
\item If $\mathfrak{q}$ is a prime in $K$ above $\mathfrak{p}$, and $\mathfrak{q}$ is inert in $\Lrho$, or if $\mathfrak{p}$ has inertia degree $3$ in $\Lrho$, then $\mathrm{Frob}_\mathfrak{p}|_{\Lrho}$ has order $3$ in $\mathrm{Gal}(\Lrho/F)$, and thus $\widetilde{\rho}(\mathrm{Frob}_\mathfrak{p})$ has odd trace.
\item If $\mathfrak{q}$ is a prime in $K$ above $\mathfrak{p}$, and $\mathfrak{q}$ splits completely in $\Lrho$, then $\mathrm{Frob}_\mathfrak{p}|_{\Lrho}$ has order at most $2$ in $\mathrm{Gal}(\Lrho/F)$, and thus $\widetilde{\rho}(\mathrm{Frob}_\mathfrak{p})$ has even trace. 
\end{itemize}

Suppose, first, that $\widetilde{\rho_E}$ has trivial image, so that $L_E = F$. If $\widetilde{\rho_\pi}$ were to have non-trivial image, then there must be an intermediate extension $K/F$ such that $L_\pi/K$ is a cubic Galois extension. From our earlier comments, we know that every possible choice for $K$ corresponds to a quadratic character of the ray class group $Cl(\mathcal{O}_F,\mathfrak{m}_{F,2})$.

Using PARI (\textbf{[PARI]}), one can determine which choice of $K$ corresponds to a given quadratic character. For each such choice, we need to show that there is no cubic Galois extension of $K$ unramified away from the primes dividing $2\mathfrak{n}$. To determine this, we construct the ray class group $Cl(\mathcal{O}_K,\mathfrak{m}_{K,3})$, and consider its cubic characters.

Each cubic extension of $K$ corresponds to such a character. Fixing a $\mz/3\mz$-basis $\{\chi_1,\ldots,\chi_n\}$ of the group of cubic characters, choose a set $\{\mathfrak{q}_1,\ldots,\mathfrak{q}_n\}$ of prime ideals in $K$ such that the vectors $(\chi_1(\mathfrak{q}_i),\ldots,\chi_n(\mathfrak{q}_i))$ span $(\mz/p\mz)^n$ (which can be shown to exist using Chebotarev's density theorem). It is straightforward to see that, for any non-trivial cubic character $\chi$ of $Cl(\mathcal{O}_K,\mathfrak{m}_{K,3})$, we must have $\chi(\mathfrak{q}_i) \neq 0$ for some prime $\mathfrak{q}_i$.

Now, if $L_\pi$ is a cubic extension of $K$, then there must be a corresponding cubic character $\chi_\pi$, and we must have $\chi_\pi(\mathfrak{q}_i) \neq 0$ for one of the primes $\mathfrak{q}_i$ chosen above. Then $\mathfrak{q}_i$ is inert in $L_\pi$, and so $\widetilde{\rho_\pi}(\mathrm{Frob}_\mathfrak{p})$ must have odd trace. 

Thus, to determine isomorphism, we compute the finite set $\{\mathfrak{q}_1,\ldots,\mathfrak{q}_n\}$ of primes as above, and find the corresponding primes $\{\mathfrak{p}_1,\ldots,\mathfrak{p}_n\}$ of $F$ lying beneath them. If $\mathrm{Tr}(\widetilde{\rho_\pi}(\mathrm{Frob}_{\mathfrak{p}_i}))$ is even for each prime $\mathfrak{p}_i$, then the extension $L_\pi/K$ must be trivial, and so $\widetilde{\rho_\pi}$ has trivial image, as required.

If $\widetilde{\rho_E}$ has $C_3$-image, we must first show that the image of $\widetilde{\rho_\pi}$ is non-trivial, and that $L_\pi$ contains no quadratic extension of $F$ (and thus must be a cubic extension of $F$), and then show that $L_\pi$ can only correspond to $L_E$.

To this end, we consider the quadratic and cubic characters of the ray class group $Cl(\mathcal{O}_F,\mathfrak{m}_{F,6})$. Firstly, let $\{\chi_1, \ldots, \chi_n\}$ be a $\mz/2\mz$-basis of the quadratic characters, and choose a set $\{\mathfrak{p}_1, \ldots, \mathfrak{p}_n\}$ of prime ideals in $F$ such that the vectors $(\chi_1(\mathfrak{p}_i),\ldots,\chi_n(\mathfrak{p}_i))$ span $(\mz/2\mz)^n$. If $L_\pi$ contained a quadratic extension $K$ of $F$, then there would be some quadratic character $\chi_K$ of $Cl(\mathcal{O}_F,\mathfrak{m}_{F,6})$ corresponding to it, and we would have $\chi_K(\mathfrak{p}_i) \neq 0$ for some prime $\mathfrak{p}_i$. In particular, $\mathfrak{p}_i$ must be inert in $K$, and so $\widetilde{\rho_\pi}(\mathrm{Frob}_{\mathfrak{p}_i})$ has even trace. If, therefore, $\mathrm{Tr}(\widetilde{\rho_\pi}(\mathrm{Frob}_{\mathfrak{p}_i}))$ is odd for each prime $\mathfrak{p}_i$, $L_\pi$ can contain no quadratic extension of $F$. In addition, since the identity in $\mathrm{GL}_2(\mf_2)$ has even trace, this implies that $L_\pi$ is a non-trivial (and hence cubic) extension of $F$.

It remains to prove that the extension $L_\pi$ is in fact $L_E$, for which we consider the cubic characters of $Cl(\mathcal{O}_F,\mathfrak{m}_{F,6})$. Identify the character $\chi_E$ corresponding to $L_E$, extend this to a $\mz/3\mz$-basis $\{\chi_E,\chi_1,\ldots,\chi_n\}$ of the cubic characters, and choose a set $\{\mathfrak{p}_1,\ldots,\mathfrak{p}_n\}$ of prime ideals in $F$ such that $\chi_E(\mathfrak{p}_i) = 0$ for each $\mathfrak{p}_i$, and the vectors $(\chi_1(\mathfrak{p}_i),\ldots,\chi_n(\mathfrak{p}_n))$ span $(\mz/3\mz)^n$. 

Now, if $L_\pi$ is not isomorphic to $L_E$, then one finds that the corresponding character $\chi_\pi$ must satisfy $\chi_\pi(\mathfrak{p}_i) \neq 0$ for some $\mathfrak{p}_i$. In this case, $\mathfrak{p}_i$ must be inert in $L_\pi$, and so $\widetilde{\rho_\pi}(\mathrm{Frob}_{\mathfrak{p}_i})$ must have odd trace. If, then, $\mathrm{Tr}(\widetilde{\rho_\pi}(\mathrm{Frob}_{\mathfrak{p}_i}))$ is even for each prime $\mathfrak{p}_i$, we deduce that the residual representations are isomorphic.

Unsurprisingly, if $\widetilde{\rho_E}$ has $S_3$-image, the task is somewhat more complicated. To begin with, note that the group $S_3$ has a unique (up to inner automorphism) subgroup of order $2$, and thus there exists an intermediate quadratic extension $K/F$ such that $L_E/K$ is a cubic Galois extension. We shall first show that if $L_\pi$ contains a quadratic extension of $F$, then it must be isomorphic to $K$, and then show that $L_\pi$ can contain no cubic extension of $F$.

To do this, we construct the ray class group $Cl(\mathcal{O}_F,\mathfrak{m}_{F,6})$, whose quadratic and cubic characters correspond to all possible quadratic and cubic extensions of $F$ contained in $L_\pi$.

Let $\chi_K$ denote the unique quadratic character of $Cl(\mathcal{O}_F,\mathfrak{m}_{F,6})$ corresponding to the field $K$, and extend this to a $\mz/2\mz$-basis $\{\chi_K,\chi_1,\ldots,\chi_n\}$ of quadratic characters of $Cl(\mathcal{O}_F,\mathfrak{m}_{F,6})$. Let $\{\mathfrak{p}_1,\ldots,\mathfrak{p}_n\}$ be prime ideals in $F$ such that the vectors $(\chi_1(\mathfrak{p}_i),\ldots,\chi_n(\mathfrak{p}_n))$ span $(\mz/2\mz)^n$ \it{and} $\chi_K(\mathfrak{p}_i) = 0$ for all $\mathfrak{p}_i$. As before, it is not too difficult to see that for \it{any} quadratic character $\chi$ which is not a multiple of $\chi_K$, we must have $\chi(\mathfrak{p}_i) \neq 0$ for some $\mathfrak{p}_i$. 

In particular, suppose $L_\pi$ contains a quadratic extension of $F$ which is \it{not} isomorphic to $K$, let $\chi$ be the corresponding quadratic character, and choose a prime $\mathfrak{p}_i$ for which $\chi(\mathfrak{p}_i) \neq 0$. Then $\mathfrak{p}_i$ must have inertia degree $2$ in $L_\pi$, and consequently must have even trace. If $\mathrm{Tr}(\widetilde{\rho_\pi}(\mathrm{Frob}_{\mathfrak{p}_i}))$ is odd for all primes $\mathfrak{p}_i$ in the finite set defined above, then we conclude that the only possible quadratic extension of $F$ contained in $L_\pi$ is $K$. We also note that, if all these traces are odd, then the image of $\widetilde{\rho_\pi}$ contains an element of order $3$.

We next need to eliminate the possibility that $\widetilde{\rho_\pi}$ has $C_3$-image. If this were the case, then we could define a cubic character $\chi_\pi$ of $Cl(\mathcal{O}_F,\mathfrak{m}_{F,6})$ corresponding to $L_\pi$. Denoting by $\{\chi_1,\ldots,\chi_n\}$ a $\mz/3\mz$-basis of cubic characters, we can find prime ideals $\{\mathfrak{p}_1,\ldots,\mathfrak{p}_n\}$ in $F$ such that the vectors $(\chi_1(\mathfrak{p}_i),\ldots,\chi_n(\mathfrak{p}_i))$ span $(\mz/3\mz)^n$ and, moreover, such that each $\mathfrak{p}_i$ either splits completely in $L_E$, or is inert in $K$. For each such prime, $\widetilde{\rho_\pi}(\mathrm{Frob}_{\mathfrak{p}_i})$ has even trace.

As usual, the non-triviality of $\chi_\pi$ implies that $\chi_\pi(\mathfrak{p}_i) \neq 0$ for some $\mathfrak{p}_i$. But this in turn implies that $\mathfrak{p}_i$ has inertia degree $3$ in $L_\pi$, and thus $\widetilde{\rho_\pi}(\mathrm{Frob}_{\mathfrak{p}_i})$ must have odd trace, so if $\mathrm{Tr}(\widetilde{\rho_\pi}(\mathrm{Frob}_{\mathfrak{p}_i}))$ is even for each prime $\mathfrak{p}_i$ above, $\widetilde{\rho_\pi}$ cannot have $C_3$-image. Since we have already deduced the existence of an element of order $3$, it must therefore have full image in $\mathrm{GL}_2(\mf_2)$.

We therefore know that $L_E$ and $L_\pi$ are both $S_3$-extensions of $F$, and both contain the quadratic extension $K$ of $F$. It remains to show that the $C_3$-extension $L_\pi/K$ is in fact $L_E$. To do this, we construct the ray class group $Cl(\mathcal{O}_K,\mathfrak{m}_{K,3})$, and consider its cubic characters.

Fix a cubic character $\chi_E$ of $Cl(\mathcal{O}_K,\mathfrak{m}_{K,3})$ corresponding to the extension $F_E$, extend this to a $\mz/3\mz$-basis $\{\chi_E,\chi_1,\ldots,\chi_n\}$ of such characters, and let $\{\mathfrak{q}_1,\ldots,\mathfrak{q}_n\}$ be prime ideals in $K$ such that the vectors $(\chi_1(\mathfrak{q}_i),\ldots,\chi_n(\mathfrak{q}_i))$ span $(\mz/3\mz)^n$, and $\chi_E(\mathfrak{q}_i) = 0$ for all $i$. Each such prime must split completely in $L_E$, and so, for any prime $\mathfrak{p}$ in $F$ lying beneath one of the $\mathfrak{q}_i$, $\widetilde{\rho_E}(\mathrm{Frob}_\mathfrak{p})$ has even trace.

If $L_\pi$ were not isomorphic to $L_E$, then we could find a cubic character $\chi_\pi$ of $Cl(\mathcal{O}_K,\mathfrak{m}_{K,3})$, corresponding to it, for which $\chi_\pi(\mathfrak{q}_i) \neq 0$ for one of the $\mathfrak{q}_i$. Consequently $\mathfrak{q}_i$ has inertia degree $3$ in $L_\pi$, and so $\widetilde{\rho_\pi}(\mathrm{Frob}_\mathfrak{p})$ must have odd trace. If, then, $\mathrm{Tr}(\widetilde{\rho_\pi}(\mathrm{Frob}_\mathfrak{p}))$ is even for each prime $\mathfrak{p}$ lying beneath one of the $\mathfrak{q}_i$, we can deduce that $L_\pi$ \it{must} be isomorphic to $L_E$, and so the residual representations are isomorphic, as required. 
\end{document}